\documentclass{article}
\usepackage{jlcode}
\RequirePackage{amsmath}
\RequirePackage{amsthm}
\RequirePackage{amsfonts}
\RequirePackage{amssymb}
\RequirePackage{mathtools}  
\RequirePackage{mathrsfs}   
\RequirePackage{bbm}        
\RequirePackage{nicefrac}   
\RequirePackage{siunitx}
\RequirePackage{url}
\RequirePackage{diagbox}
\usepackage{arxiv}
\usepackage{kthcolors}
\usepackage{algorithm}
\usepackage{algpseudocode}

\algnewcommand\algorithmicinput{\hspace*{1.0em}\textbf{input:}}
\algnewcommand\Input{\item[\algorithmicinput]}
\algnewcommand\algorithmicoutput{\hspace*{1.0em}\textbf{output:}}
\algnewcommand\Output{\item[\algorithmicoutput]}

\DeclarePairedDelimiter{\parentheses}{(}{)}
\DeclarePairedDelimiter{\brackets}{[}{]}
\DeclarePairedDelimiter{\braces}{\lbrace}{\rbrace}
\DeclarePairedDelimiter{\verts}{\lvert}{\rvert}

\DeclarePairedDelimiterX{\Set}[1]{\lbrace}{\rbrace}{
  
  #1
}

\providecommand{\st}{\ensuremath\;\mathrm{s.t.}\;}
\providecommand{\expect}[2][{}k]{\ensuremath{\operatorname{\mathbb{E}}_{#1}\brackets*{#2}}}
\providecommand{\abs}[1]{\ensuremath{\verts*{#1}}}

\providecommand{\argmin}[2][x]{\ensuremath{\arg\,\min_{#1}\braces*{#2}}}

\DeclareMathOperator*{\minimize}{minimize}
\DeclareMathOperator*{\maximize}{maximize}

\theoremstyle{definition}                     

\theoremstyle{remark}                         

\theoremstyle{plain}                          

\usepackage{pgfplots}
\usepackage{tikz}
\usetikzlibrary{arrows.meta,intersections,fit,math,shapes,calc,positioning}

\begin{document}

\title{Optimal day-ahead orders using stochastic programming and noise-driven RNNs}

\author{Martin Biel \\
  Division of Decision and Control Systems\\
  School of EECS, KTH Royal Institute of Technology\\
  SE-100 44 Stockholm, Sweden\\
  \texttt{mbiel@kth.se}}

\date{\today}

\maketitle

\begin{abstract}
  This paper presents a methodology for strategic day-ahead planning that uses a combination of deep learning and optimization. A noise-driven recurrent neural network structure is proposed for forecasting electricity prices and local inflow to water reservoirs. The resulting forecasters generate predictions with seasonal variation without relying on long input sequences. This forecasting method is employed in a stochastic program formulation of the day-ahead problem. This results in optimal order strategies for a price-taking hydropower producer participating in the Nordic day-ahead market. Using an open-source software framework for stochastic programming, the model is implemented and distributed over multiple cores. The model is then solved in parallel using a sampling-based algorithm. Tight confidence intervals around the stochastic solution are provided, which show that the gain from adopting a stochastic approach is statistically significant.
\end{abstract}

\section{Introduction}
\label{sec:ch5-introduction}

Electricity trading in the Nordic energy market is primarily based on day-ahead auctions. Producers and consumers submit orders of price and electricity volumes to convey their desired market participation for the upcoming day. The market operator then settles a market price based on supply and demand in the submitted orders. Since prices and demands are unknown when submitting orders, imbalances can occur. These can then be mitigated in various short-term markets. Hydropower producers participating in these markets have the opportunity to store energy in the reservoirs, which leads to high flexibility. Due to the high uncertainty associated with day-ahead market participation, precise forecasts and planning procedures are required to devise optimal order strategies. Stochastic programming is an effective approach for modeling these types of problems where decisions subjected to uncertainty are taken in stages.

We consider a complete process where we forecast, model, and solve, in order to determine optimal order strategies for a hydropower producer participating in the Nordic day-ahead market NordPool. We formulate the day-ahead planning problem as a stochastic program and implement it in the open-source software framework, \jlinl{StochasticPrograms.jl}~\cite{spjl}, written in the Julia programming language. The framework is implemented with parallel capabilities, which allows us to efficiently distribute stochastic programs in memory. We use these features to instantiate large-scale day-ahead models. We use historical data to construct samplers, which in turn generate the scenarios used to construct the stochastic programs. The samplers are based on the proposed noise-driven recurrent neural network (RNN) structure, which is used to forecast day-ahead prices and local inflows to water reservoirs. A key feature is that the samplers generate forecasts without relying on long input sequences. Also, they can account for seasonal variations in the forecasts through exogenous input parameters. To generate stochastic solutions of the formulated day-ahead model with tight confidence intervals, we apply a \emph{sampled average approximation} (SAA) algorithm. This allows us to verify that the \emph{value of the stochastic solution} (VSS) of the model is statistically significant in relation to the model and forecasts. We employ a family of distributed L-shaped algorithms~\cite{distlshaped} to efficiently run SAA in parallel on a computational cluster. We make use of acceleration techniques to improve the convergence of the algorithm.

The application of stochastic programming for day-ahead planning has been actively studied. Some notable papers are~\cite{prosumer, Fleten2007, Asgard2018}. Also, an informative survey is provided by the authors of~\cite{Asgard2019}. A common shortcoming in these contributions is that the resulting VSS is relatively low. Moreover, it is common that low VSS numbers are presented without acknowledging that they might be within the error bounds of the stochastic solution. In this work, we carefully ensure that the presented VSS is statistically significant in relation to the model we solve.

RNN type approaches have been successfully applied to both price forecasting~\cite{rnnprice, rnnprice2} and inflow forecasting~\cite{rnninflow1,rnninflow2}. A common characteristic in these implementations is that the forecasts rely on a lagged input sequence of historical values. Consequently, the network training involves long input sequences and performance can be impacted negatively if another historical sequence is given as input after training. This reduces the generality of the forecasts. This also holds for ARIMA and ARMAX type models~\cite{priceforecasting}. The proposed noise-driven RNN forecaster does not have this characteristic. The networks are trained on historical data, but forecasts are generated from Gaussian inputs.

\section{The day-ahead problem}
\label{sec:day-ahead-problem}

A day-ahead planning problem involves specifying optimal order volumes in a deregulated electricity market. We give a brief introduction to this problem here, as well as narrow down the scope of the problem solved and presented in this contribution.

\subsection{The day-ahead market}
\label{sec:day-ahead-market}

In a deregulated day-ahead electricity market, producers and consumers place orders that specify the electricity volumes they wish to sell and buy the next day. The next-day market price is determined by the equilibrium price of these orders. After market price settlement, all volume orders that have a price equal to or lower than then market price are accepted. All market participants then become balance responsible for their accepted orders. Any imbalances can be continuously adjusted by submitting orders to an intraday market after the market price settlement. In general, imbalance settlement will involve a less favorable price than the day-ahead price. Since the next-day market price is unknown, and the cost of imbalances can be high, careful planning is required in order to submit strategic orders to the day-ahead market.

\subsection{Order types}
\label{sec:order-types}

The Nordic day-ahead market offers four order variants for trading electricity volumes, hourly orders, block orders, exclusive groups, and flexible orders. We give a brief introduction to hourly orders and regular block orders.

\begin{figure}
  \centering
  \resizebox{\textwidth}{!}{
\begin{tikzpicture}[]
\begin{axis}[height = {101.6mm}, legend pos = {north west}, ylabel = {Price [EUR/MWh]}, title = {Order Curve}, xmin = {-48.57420166709518}, xmax = {728.6130250064277}, ymax = {34.646125008448806}, xlabel = {Order Volume [MWh/h]}, unbounded coords=jump,scaled x ticks = false,xlabel style = {font = {\fontsize{11 pt}{14.3 pt}\selectfont}, color = {rgb,1:red,0.00000000;green,0.00000000;blue,0.00000000}, draw opacity = 1.0, rotate = 0.0},xmajorgrids = true,xtick = {0.0,113.33980388988876,226.6796077797775,340.0194116696663,453.359215559555,566.6990194494439,680.0388233393326},xticklabels = {0.00,113.34,226.68,340.02,453.36,566.70,680.04},xtick align = inside,xticklabel style = {font = {\fontsize{8 pt}{10.4 pt}\selectfont}, color = {rgb,1:red,0.00000000;green,0.00000000;blue,0.00000000}, draw opacity = 1.0, rotate = 0.0},x grid style = {color = {rgb,1:red,0.00000000;green,0.00000000;blue,0.00000000},
draw opacity = 0.1,
line width = 0.5,
solid},axis x line* = left,x axis line style = {color = {rgb,1:red,0.00000000;green,0.00000000;blue,0.00000000},
draw opacity = 1.0,
line width = 1,
solid},scaled y ticks = false,ylabel style = {font = {\fontsize{11 pt}{14.3 pt}\selectfont}, color = {rgb,1:red,0.00000000;green,0.00000000;blue,0.00000000}, draw opacity = 1.0, rotate = 0.0},ymajorgrids = true,ytick = {0.0,2.7505129689722754,5.501025937944551,8.251538906916826,11.002051875889101,13.752564844861377,16.503077813833652,19.25359078280593,22.004103751778203,24.754616720750477,27.505129689722754,30.25564265869503},yticklabels = {0.00,2.75,5.50,8.25,11.00,13.75,16.50,19.25,22.00,24.75,27.51,30.26},ytick align = inside,yticklabel style = {font = {\fontsize{8 pt}{10.4 pt}\selectfont}, color = {rgb,1:red,0.00000000;green,0.00000000;blue,0.00000000}, draw opacity = 1.0, rotate = 0.0},y grid style = {color = {rgb,1:red,0.00000000;green,0.00000000;blue,0.00000000},
draw opacity = 0.1,
line width = 0.5,
solid},axis y line* = left,y axis line style = {color = {rgb,1:red,0.00000000;green,0.00000000;blue,0.00000000},
draw opacity = 1.0,
line width = 1,
solid},    xshift = 0.0mm,
    yshift = 0.0mm,
    axis background/.style={fill={rgb,1:red,1.00000000;green,1.00000000;blue,1.00000000}}
,title style = {font = {\fontsize{14 pt}{18.2 pt}\selectfont}, color = {rgb,1:red,0.00000000;green,0.00000000;blue,0.00000000}, draw opacity = 1.0, rotate = 0.0},legend style = {color = {rgb,1:red,0.00000000;green,0.00000000;blue,0.00000000},
draw opacity = 1.0,
line width = 1,
solid,fill = {rgb,1:red,1.00000000;green,1.00000000;blue,1.00000000},font = {\fontsize{8 pt}{10.4 pt}\selectfont}},colorbar style={title=}, ymin = {-2.7505129689722754}, width = {152.4mm}]\addplot+ [color = {rgb,1:red,0.00000000;green,0.00000000;blue,0.00000000},
draw opacity = 1.0,
line width = 2,
solid,mark = none,
mark size = 2.0,
mark options = {
    color = {rgb,1:red,0.00000000;green,0.00000000;blue,0.00000000}, draw opacity = 1.0,
    fill = {rgb,1:red,0.00000000;green,0.60560316;blue,0.97868012}, fill opacity = 1.0,
    line width = 1,
    rotate = 0,
    solid
},forget plot]coordinates {
(0.0, 0.0)
(0.0, 0.0)
(0.0, 20.89356016358743)
(636.7057290307187, 20.89356016358743)
(636.7057290307187, 20.89356016358743)
(636.7057290307187, 23.644073132559704)
(636.7057290307187, 23.644073132559704)
(636.7057290307187, 26.39458610153198)
(636.7057290307187, 29.14509907050426)
(680.0388233393326, 29.14509907050426)
(680.0388233393326, 29.14509907050426)
(680.0388233393326, 31.895612039476532)
(680.0388233393326, 31.895612039476532)
(680.0388233393326, 34.646125008448806)
};
\addplot+[draw=none, color = {rgb,1:red,0.88887350;green,0.43564919;blue,0.27812294},
draw opacity = 1.0,
line width = 0,
solid,mark = *,
mark size = 2.0,
mark options = {
    color = {rgb,1:red,0.00000000;green,0.00000000;blue,0.00000000}, draw opacity = 1.0,
    fill = {rgb,1:red,0.38431373;green,0.57254902;blue,0.18039216}, fill opacity = 1.0,
    line width = 1,
    rotate = 0,
    solid
}] coordinates {
(0.0, 0.0)
};
\addlegendentry{Price Independent Order}
\addplot+[draw=none, color = {rgb,1:red,0.24222430;green,0.64327509;blue,0.30444865},
draw opacity = 1.0,
line width = 0,
solid,mark = *,
mark size = 2.0,
mark options = {
    color = {rgb,1:red,0.00000000;green,0.00000000;blue,0.00000000}, draw opacity = 1.0,
    fill = {rgb,1:red,0.09803922;green,0.32941176;blue,0.65098039}, fill opacity = 1.0,
    line width = 1,
    rotate = 0,
    solid
}] coordinates {
(636.7057290307187, 20.89356016358743)
(636.7057290307187, 23.644073132559704)
(636.7057290307187, 26.39458610153198)
(680.0388233393326, 29.14509907050426)
(680.0388233393326, 31.895612039476532)
};
\addlegendentry{Price Dependent Order}
\addplot+ [color = {rgb,1:red,0.00000000;green,0.00000000;blue,0.00000000},
draw opacity = 1.0,
line width = 1,
dashed,mark = none,
mark size = 2.0,
mark options = {
    color = {rgb,1:red,0.00000000;green,0.00000000;blue,0.00000000}, draw opacity = 1.0,
    fill = {rgb,1:red,0.76444018;green,0.44411178;blue,0.82429754}, fill opacity = 1.0,
    line width = 1,
    rotate = 0,
    solid
},forget plot]coordinates {
(636.7057290307187, 26.39458610153198)
(680.0388233393326, 29.14509907050426)
};
\addplot+[draw=none, color = {rgb,1:red,0.67554396;green,0.55566233;blue,0.09423434},
draw opacity = 1.0,
line width = 0,
solid,mark = diamond*,
mark size = 2.0,
mark options = {
    color = {rgb,1:red,0.00000000;green,0.00000000;blue,0.00000000}, draw opacity = 1.0,
    fill = {rgb,1:red,0.38431373;green,0.57254902;blue,0.18039216}, fill opacity = 1.0,
    line width = 1,
    rotate = 0,
    solid
}] coordinates {
(661.9983026893217, 28.0)
};
\addlegendentry{Trading Outcome}
\end{axis}
\end{tikzpicture}}
  \caption{Single hourly order example, showing volume interpolation after market price settlement. The price independent order is always accepted at market price.}
  \label{fig:hourlyorderex}
\end{figure}
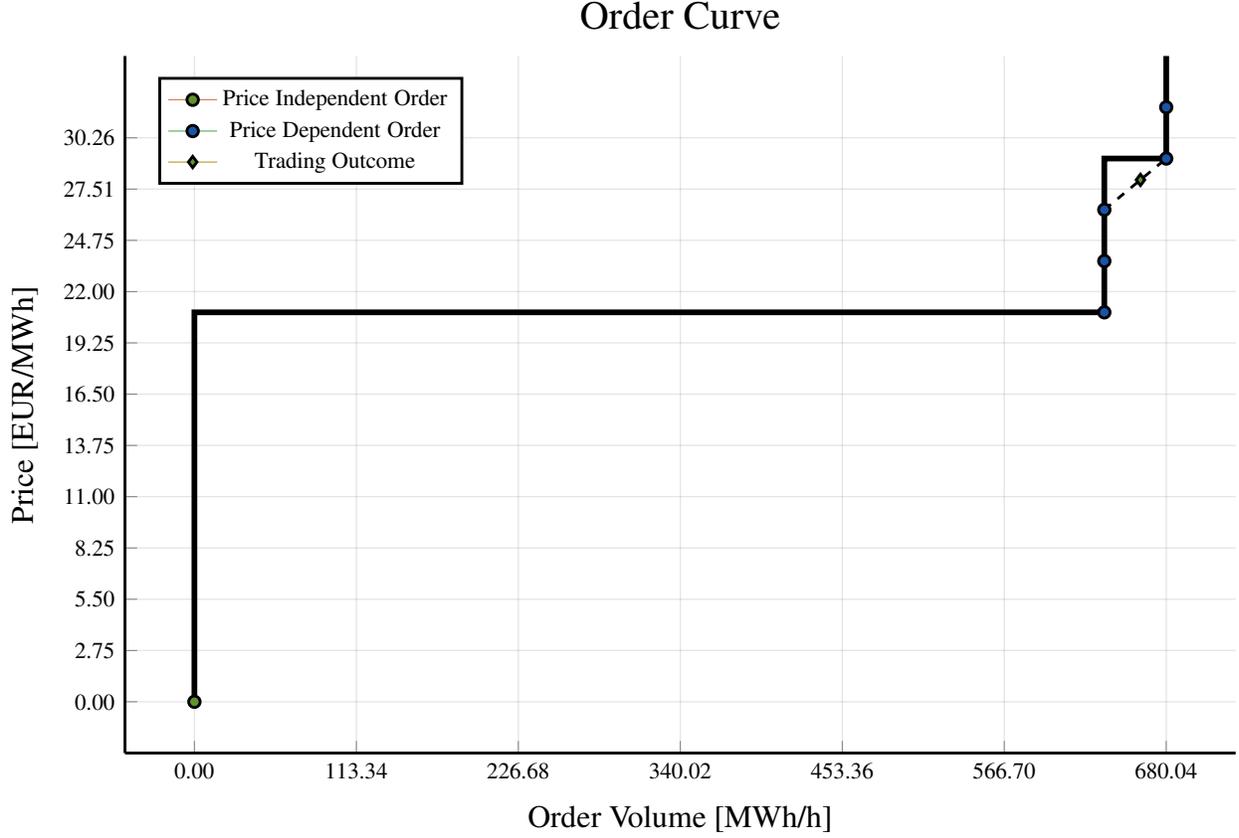

Hourly orders can be placed in two ways. A price independent hourly order specifies an electricity volume that is to be purchased or sold at market price during a certain hour, independent of the market price. A price dependent order specifies electricity volumes at given price points. If the settled market price ends up between specified price steps a linear interpolation is performed between the adjacent volume orders to determine the order volume. A settled hourly order is illustrated in Fig.~\ref{fig:hourlyorderex}.

Block orders span over an interval of consecutive hours. A regular block order is accepted in its entirety if the mean market price in the specified interval is higher or equal to the order price. The participants then become balance responsible for the order volume every hour of the specified interval, at the mean market price in the interval. Block orders where the price is higher than the mean market price in the given interval are rejected. Other block order variants exist, such as linked block orders and curtailable block orders. These all include conditional elements, and determining optimal orders would involve combinatorial optimization. This is also true for the remaining order types: exclusive group and flexible orders. We do not give further details into the conditional order types as only hourly orders and regular block orders will be used subsequently.

\subsection{Problem setting}
\label{sec:problem-setting}

In this work, we formulate a day-ahead planning problem from the perspective of a hydropower producer that is assumed to own all 15 power stations in the Swedish river Skellefteälven. The reservoirs in Skellefteälven are of varying size, ranging from smaller reservoirs by the sea outlet to complete lakes at the basin of the river. This enables flexible production schedules but requires complex planning. The producer is also assumed to be price-taking so submitted orders do not influence the market price. The day-ahead model is limited to include only hourly orders and regular block orders. When the market price has been realized, the producer optimizes the power production with respect to the price and the opportunity cost of discharging water from the reservoirs. We assume that there are no fixed contracts to adhere to. In other words, all electricity production is sold for profit in the market. All submitted orders must follow the trade regulations specified by the Nordpool market. A general description of the day-ahead problem is given in~\eqref{eq:generaldayahead}.

\begin{equation} \label{eq:generaldayahead}
\begin{aligned}
 \maximize_{\mathclap{\text{Order strategy}}} & \quad \text{Profit} + \text{Water value} - \text{Imbalances} \\
 \text{subject to} & \quad \text{Physical limitations} \\
 & \quad \text{Economic/legal limitations}
\end{aligned}
\end{equation}

Because next-day market prices are unknown when placing orders, we formulate a two-stage stochastic program to generate optimal orders. The first-stage decisions are the orders submitted to the day-ahead market. A general description of the first-stage problem is given in~\eqref{eq:generalfirststage}.

\begin{equation} \label{eq:generalfirststage}
\begin{aligned}
 \maximize_{\mathclap{\text{Order strategy}}} & \quad \expect[]{\text{Revenue}\parentheses*{\text{Order strategy}, \text{Price}, \text{Inflow}}} \\
 \text{subject to} & \quad \text{Trade regulations} \\
\end{aligned}
\end{equation}

In each second-stage scenario, uncertain parameters are realized and the electricity production is optimized with respect to profits and water value while satisfying the settled order commitments. A general description of the second stage is given in~\eqref{eq:generalsecondstage}
\begin{equation} \label{eq:generalsecondstage}
\begin{aligned}
 \maximize_{\mathclap{\text{Production schedule}}} & \quad \text{Profit}\parentheses*{\text{Price}} + \text{Water value} - \text{Imbalances} \\
 \text{subject to} & \quad \text{Commitments}\parentheses*{\text{Order strategy}, \text{Price}} \\
 & \quad \text{Hydrological balance}\parentheses*{\text{Inflow}} \\
 & \quad \text{Electricity production} \\
 & \quad \text{Load balance}\parentheses*{\text{Commitments}}
\end{aligned}
\end{equation}

In addition, the producer can take recourse decisions by trading surplus or shortage in a balancing market. Both market prices and local water inflows to the reservoirs are considered uncertain. Apart from short-term randomness, both parameters are subjected to large seasonal variations. Electricity prices in the Nordic region are generally higher in the winter period and the local inflows vary considerably over the year due to snowmelt.

\subsection{Problem parameters}
\label{sec:problem-parameters}

We now list the set of parameters required to define the day-ahead problem. We introduce deterministic parameters and uncertain parameters separately.

The deterministic parameters constitute physical hydro plant parameters and trade regulations. Physical parameters for the power stations in Skellefteälven are available in~\cite{Sag2018}. These include river topology, reservoir capacities, discharge limits, and water travel time between adjacent stations. Trade regulations, such as trading fees and order limits, are available at NordPool~\cite{nordpool}.

The uncertain parameters consist of market prices and local inflows to the reservoirs. Historical data is available for both quantities and can be used to create statistical models for the random parameters. Historical price data between 2013 and 2018 is available on NordPool~\cite{nordpooldata}. Historical local inflow data in Skellefteälven between 1999 and 2018 was acquired from the Swedish Meteorological and Hydrological Institute (SMHI), using the S-HYPE model~\cite{shype}. The model provides daily hydro flow measurement readings at every river location by a nearby observation node.

\section{Implementation}
\label{sec:impl-deta}

In this section, we go through the implementation details of this work.

\subsection{Noise-driven RNN forecasting}
\label{sec:rnn-design}

We give a general overview of the noise-driven RNN design used to create the price and inflow forecasters. As is common in deep learning applications, hyperparameter tuning involves a lot of trial and error. Therefore, we only present the final network architectures and the relevant hyperparameter values. All networks are implemented and trained using Flux.jl~\cite{Flux}.

\subsubsection{Network structure}
\label{sec:network-structure}

\begin{figure}
  \centering
  \begin{tikzpicture}
  \node at (0,0) (anchor) {};

  \node[left of=anchor,xshift=-100pt] (gaussian) {$\tilde{w}$};
  \node[draw,right of = gaussian, node distance = 60pt] (init) {Initializer network};
  \draw[->] (gaussian) -- (init);
  \node[below of = init, node distance = 20pt] (seasonal) {$s$};
  \draw[->] (seasonal) -- (init);
  \node[draw, right of = init, node distance = 100pt] (seq) {Sequence generator};
  \draw[->] (init) -- node[pos=0.5, above] {$x_0$} (seq);
  \node[below of = seq, node distance = 20pt] (seasonalseq) {$s$};
  \draw[->] (seasonalseq) -- (seq);
  \node[right of = seq, node distance = 60pt] (output) {$x_t$};
  \draw[->] (seq) -- node[pos = 0.5] (sr0) {} (output);
  \node[above of = sr0] (sr1) {};
  \node[above of = seq] (sr2) {};
  \draw[->] (sr0.center) -- (sr1.south) -- (sr2.south) -- (seq);
\end{tikzpicture}
  \caption{Noise-driven RNN forecaster architecture.}
  \label{fig:rnnarch}
\end{figure}
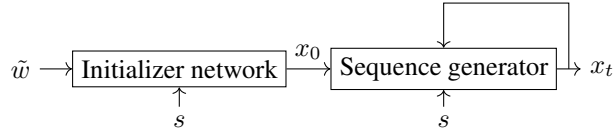

The price forecaster and the inflow forecaster share a similar structure. The main aim of this proposed structure is to enable forecasting of sequential data with seasonal variation, without having to rely on long input sequences. The general structure of the forecaster is shown in Figure~\ref{fig:rnnarch}.

The proposed network structure consists of two key components. First, an initializer network is used to compute the initial state of the forecasted sequence. The inputs to the initializer network are a set of seasonal indicators and a Gaussian noise signal. The initializer has a deep structure with a set of dense layers, all having the standard form
\begin{equation*}
  y = \sigma (Wx + b).
\end{equation*}
\noindent
Some layers use rectified linear units as their activation functions to introduce nonlinearity to the network.

Next, the computed initial state is passed to a sequence generator, which is a network that includes a recurrent layer of the form
\begin{equation*}
  \begin{aligned}
    h_t &= f(x_t, h_{t-1}) \\
    y_t &= g(h_t).
  \end{aligned}
\end{equation*}
Here, $x_t$ is the input, $h_t$ is a hidden persistent state, and $y_t$ is the output. The memory state $h_t$ enables sequence learning. We construct the recurrent layer using Gated Recurrent Units (GRU), first introduced in~\cite{gru}. A gated unit includes an \emph{update} state that decides what parts of the memory to keep, and a \emph{forget} state that decides what parts of the memory to drop. The use of GRU based networks has proven useful in similar applications~\cite{rnnprice}. The inputs to the sequence generator are given by the previous sequence value and seasonal indicators. Here, the seasonal indicators also include the current timestep of the sequence. This is used to increase the chance of learning sequential patterns in the data. Similar to the initializer, the sequence generator also includes a set of dense layers. In both the initializer network and the sequence generator, a dropout mechanism is added between layers. Consequently, randomness is present in the network through both the Gaussian input to the initializer and the random removal of neural links. This makes training more difficult but increases the generality of the resulting network.

In brief, a forecasted sequence of length $T$ is created in two steps. The initial value $x_0$ is computed by the initializer network. Afterward, the remaining $T-1$ values are computed by letting data flow recursively through the sequence generator $T-1$ times, starting with $x_0$. This is the basis for both the price forecaster and the inflow forecaster.

\subsubsection{Network architectures}
\label{sec:netw-arch}

\begin{figure}
  \centering
  \begin{tikzpicture}
  \node at (0,0) (anchor) {};

  \node[left of=anchor,xshift=-100pt] (init) {};
  \node[above of=init, node distance = 10pt] (gaussian) {$\tilde{w}$};
  \node[below of=init, node distance = 10pt] (month) {$m$};
  \node[right of = gaussian, node distance = 60pt] (l1top) {};
  \node[right of = month, node distance = 60pt] (l1bottom) {};
  \node[right of = init, node distance = 60pt] (l1) {$2 \to 128$};
  \node[draw,fit=(l1top)(l1bottom)(l1)] (l1box) {};
  \node[above of = l1top, node distance = 15pt] {Dense};
  \draw[->] (gaussian) -- ($(l1box) + (-25pt,10pt)$);
  \draw[->] (month) -- ($(l1box) + (-25pt,-10pt)$);
  \node[right of = l1top, node distance = 60pt] (l2top) {};
  \node[right of = l1bottom, node distance = 60pt] (l2bottom) {};
  \node[right of = l1, node distance = 60pt] (l2) {$128 \to 128$};
  \node[draw, fit=(l2top)(l2bottom)(l2)] (l2box) {};
  \node[above of = l2top, node distance = 15pt] {Dense};
  \draw[->,dashed] (l1box) -- (l2box);
  \begin{scope}[transform canvas={yshift=7pt}]
    \draw[->,dashed] (l1box) -- (l2box);
  \end{scope}
  \begin{scope}[transform canvas={yshift=15pt}]
    \draw[->,dashed] (l1box) -- node[pos=0.5, above] (d1) {$0.5$} (l2box);
  \end{scope}
  \begin{scope}[transform canvas={yshift=-7pt}]
    \draw[->,dashed] (l1box) -- (l2box);
  \end{scope}
  \begin{scope}[transform canvas={yshift=-15pt}]
    \draw[->,dashed] (l1box) -- (l2box);
  \end{scope}
  \node[right of = l2top, node distance = 65pt] (l3top) {};
  \node[right of = l2bottom, node distance = 65pt] (l3bottom) {};
  \node[right of = l2, node distance = 65pt] (l3) {$128 \to 128$};
  \node[above of = l3top, node distance = 15pt] {Dense};
  \node[draw, fit=(l3top)(l3bottom)(l3)] (l3box) {};
  \draw[->,dashed] (l2box) -- (l3box);
  \begin{scope}[transform canvas={yshift=7pt}]
    \draw[->,dashed] (l2box) -- (l3box);
  \end{scope}
  \begin{scope}[transform canvas={yshift=15pt}]
    \draw[->,dashed] (l2box) -- node[pos=0.5, above] (d2) {$0.5$} (l3box);
  \end{scope}
  \node[above of = d2, node distance = 25pt] {relu};
  \begin{scope}[transform canvas={yshift=-7pt}]
    \draw[->,dashed] (l2box) -- (l3box);
  \end{scope}
  \begin{scope}[transform canvas={yshift=-15pt}]
    \draw[->,dashed] (l2box) -- (l3box);
  \end{scope}
  \node[right of = l3top, node distance = 60pt] (l4top) {};
  \node[right of = l3bottom, node distance = 60pt] (l4bottom) {};
  \node[right of = l3, node distance = 60pt] (l4) {$128 \to 1$};
  \node[above of = l4top, node distance = 15pt] {Dense};
  \node[draw, fit=(l4top)(l4bottom)(l4)] (l4box) {};
  \draw[->,dashed] (l3box) -- (l4box);
  \begin{scope}[transform canvas={yshift=7pt}]
    \draw[->,dashed] (l3box) -- (l4box);
  \end{scope}
  \begin{scope}[transform canvas={yshift=15pt}]
    \draw[->,dashed] (l3box) -- node[pos=0.5, above] (d3) {$0.5$} (l4box);
  \end{scope}
  \node[above of = d3, node distance = 25pt] {relu};
  \begin{scope}[transform canvas={yshift=-7pt}]
    \draw[->,dashed] (l3box) -- (l4box);
  \end{scope}
  \begin{scope}[transform canvas={yshift=-10pt}]
    \draw[->,dashed] (l3box) -- (l4box);
  \end{scope}
  \node[right of=l4, node distance = 60pt] (output) {$x_0$};
  \draw[->] (l4box) -- (output);
\end{tikzpicture}
  \caption{Initializer in the price forecaster. Activation functions and dropout values between layers are given above the neuron links.}
  \label{fig:priceinitarch}
\end{figure}
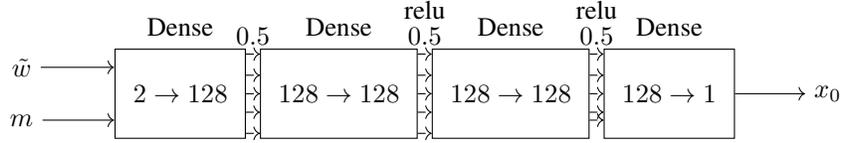

The price forecaster outputs a $24$-hour price sequence for the upcoming day. The initializer architecture is shown in Figure~\ref{fig:priceinitarch}. The output layer combines the values of the neurons in the final hidden layer into a proposed initial value. Apart from the one-dimensional Gaussian input, the initializer is also supplied with the current month.

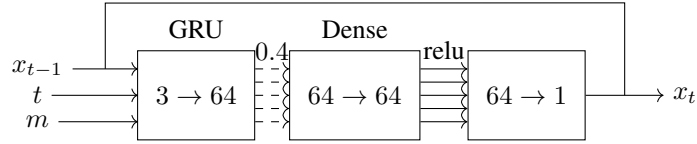
\begin{figure}
  \centering
  \begin{tikzpicture}
  \node at (0,0) (anchor) {};

  \node[left of = anchor,xshift=-100pt] (hour) {$t$};
  \node[above of = hour, node distance = 10pt] (state) {$x_{t-1}$};
  \node[below of = hour, node distance = 10pt] (month) {$m$};
  \node[right of = state, node distance = 60pt] (l1top) {};
  \node[right of = month, node distance = 60pt] (l1bottom) {};
  \node[right of = hour, node distance = 60pt] (l1) {$3 \to 64$};
  \node[draw,fit=(l1top)(l1bottom)(l1)] (l1box) {};
  \node[above of = l1top, node distance = 15pt] {GRU};
  \draw[->] (state) -- node[pos = 0.5] (sre) {} ($(l1box) + (-22pt,10pt)$);
  \draw[->] (hour) -- ($(l1box) + (-22pt,0)$);
  \draw[->] (month) -- ($(l1box) + (-22pt,-10pt)$);
  \node[right of = l1top, node distance = 60pt] (l2top) {};
  \node[right of = l1bottom, node distance = 60pt] (l2bottom) {};
  \node[right of = l1, node distance = 60pt] (l2) {$64 \to 64$};
  \node[draw, fit=(l2top)(l2bottom)(l2)] (l2box) {};
  \node[above of = l2top, node distance = 15pt] {Dense};
  \draw[->,dashed] (l1box) -- (l2box);
  \begin{scope}[transform canvas={yshift=5pt}]
    \draw[->,dashed] (l1box) -- (l2box);
  \end{scope}
  \begin{scope}[transform canvas={yshift=10pt}]
    \draw[->,dashed] (l1box) -- node[pos=0.5, above] {$0.4$} (l2box);
  \end{scope}
  \begin{scope}[transform canvas={yshift=-5pt}]
    \draw[->,dashed] (l1box) -- (l2box);
  \end{scope}
  \begin{scope}[transform canvas={yshift=-10pt}]
    \draw[->,dashed] (l1box) -- (l2box);
  \end{scope}
  \node[right of = l2top, node distance = 65pt] (l3top) {};
  \node[right of = l2bottom, node distance = 65pt] (l3bottom) {};
  \node[right of = l2, node distance = 65pt] (l3) {$64 \to 1$};
  \node[draw, fit=(l3top)(l3bottom)(l3)] (l3box) {};
  \draw[->] (l2box) -- (l3box);
  \begin{scope}[transform canvas={yshift=5pt}]
    \draw[->] (l2box) -- (l3box);
  \end{scope}
  \begin{scope}[transform canvas={yshift=10pt}]
    \draw[->] (l2box) -- node[pos=0.5, above] {relu} (l3box);
  \end{scope}
  \begin{scope}[transform canvas={yshift=-5pt}]
    \draw[->] (l2box) -- (l3box);
  \end{scope}
  \begin{scope}[transform canvas={yshift=-10pt}]
    \draw[->] (l2box) -- (l3box);
  \end{scope}
  \node[right of=l3, node distance = 60pt] (output) {$x_t$};
  \draw[->] (l3box) -- node[pos = 0.5] (sr0) {} (output);
  \node[above of = sr0, node distance = 35pt] (sr1) {};
  \node[above of = sre] (sr2) {};
  \draw (sr0.center) -- (sr1.center) -- (sr2.south) -- (sre.center);
\end{tikzpicture}
  \caption{Sequence generation network for the price forecaster. The output is fed back into the system until the full sequence is generated. The seasonal parameters are given by the current month and hour.}
  \label{fig:priceseqarch}
\end{figure}

The sequence generator is a three-layered network, as shown in Fig~\ref{fig:priceseqarch}. The first layer is a GRU based layer with $64$ hidden neurons. It is followed by a dense layer of $64$ hidden neurons and ReLU activation functions. A dropout mechanism of value $0.4$ is applied between the recurrent layer and the dense layer. The output layer combines the neuron values of the dense layer into the next sequence value, which is also fed back into the network until the full sequence has been computed. Apart from the previous sequence value, the sequence generator is also supplied with the current month and the current hour. A code excerpt exemplifying how the components of the price forecaster is created in Flux.jl is given in Listing~\ref{lst:fluxdef}.

\begin{lstlisting}[language = julia, caption = {Definition of the price forecaster components in Flux.jl}, label = {lst:fluxdef}]
initializer = Chain(Dense(2,128),
                    Dropout(0.5),
                    Dense(128,128,relu),
                    Dropout(0.5),
                    Dense(128,128,relu),
                    Dropout(0.5),
                    Dense(128, 1))

sequence_generator = Chain(GRU(3, 64),
                           Dropout(0.4),
                           Dense(64,64,relu),
                           Dense(64, 1))
\end{lstlisting}The sequence generator is a three-layered network, as shown in Fig~\ref{fig:priceseqarch}. The first layer is a GRU based layer with $64$ hidden neurons. It is followed by a dense layer of $64$ hidden neurons and ReLU activation functions. A dropout mechanism of value $0.4$ is applied between the recurrent layer and the dense layer. The output layer combines the neuron values of the dense layer into the next sequence value, which is also fed back into the network until the full sequence has been computed. Apart from the previous sequence value, the sequence generator is also supplied with the current month and the current hour. A code excerpt exemplifying how the components of the price forecaster is created in Flux.jl is given in Listing~\ref{lst:fluxdef}.

\begin{lstlisting}[language = julia, caption = {Definition of the price forecaster components in Flux.jl}, label = {lst:fluxdef}]
initializer = Chain(Dense(2,128),
                    Dropout(0.5),
                    Dense(128,128,relu),
                    Dropout(0.5),
                    Dense(128,128,relu),
                    Dropout(0.5),
                    Dense(128, 1))

sequence_generator = Chain(GRU(3, 64),
                           Dropout(0.4),
                           Dense(64,64,relu),
                           Dense(64, 1))
\end{lstlisting}

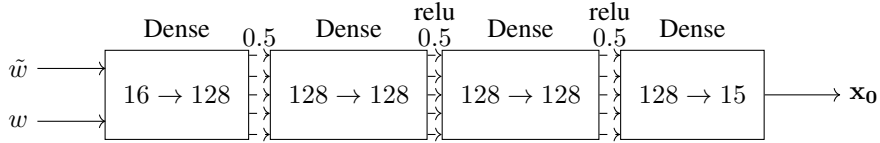
\begin{figure}
  \centering
  \begin{tikzpicture}
  \node at (0,0) (anchor) {};

  \node[left of=anchor,xshift=-100pt] (init) {};
  \node[above of=init, node distance = 10pt] (gaussian) {$\tilde{w}$};
  \node[below of=init, node distance = 10pt] (week) {$w$};
  \node[right of = gaussian, node distance = 60pt] (l1top) {};
  \node[right of = week, node distance = 60pt] (l1bottom) {};
  \node[right of = init, node distance = 60pt] (l1) {$16 \to 128$};
  \node[draw,fit=(l1top)(l1bottom)(l1)] (l1box) {};
  \node[above of = l1top, node distance = 15pt] {Dense};
  \draw[->] (gaussian) -- ($(l1box) + (-27pt,10pt)$);
  \draw[->] (week) -- ($(l1box) + (-27pt,-10pt)$);
  \node[right of = l1top, node distance = 65pt] (l2top) {};
  \node[right of = l1bottom, node distance = 65pt] (l2bottom) {};
  \node[right of = l1, node distance = 65pt] (l2) {$128 \to 128$};
  \node[draw, fit=(l2top)(l2bottom)(l2)] (l2box) {};
  \node[above of = l2top, node distance = 15pt] {Dense};
  \draw[->,dashed] (l1box) -- (l2box);
  \begin{scope}[transform canvas={yshift=7pt}]
    \draw[->,dashed] (l1box) -- (l2box);
  \end{scope}
  \begin{scope}[transform canvas={yshift=15pt}]
    \draw[->,dashed] (l1box) -- node[pos=0.5, above] (d2) {$0.5$} (l2box);
  \end{scope}
  \begin{scope}[transform canvas={yshift=-7pt}]
    \draw[->,dashed] (l1box) -- (l2box);
  \end{scope}
  \begin{scope}[transform canvas={yshift=-15pt}]
    \draw[->,dashed] (l1box) -- (l2box);
  \end{scope}
  \node[right of = l2top, node distance = 65pt] (l3top) {};
  \node[right of = l2bottom, node distance = 65pt] (l3bottom) {};
  \node[right of = l2, node distance = 65pt] (l3) {$128 \to 128$};
  \node[above of = l3top, node distance = 15pt] {Dense};
  \node[draw, fit=(l3top)(l3bottom)(l3)] (l3box) {};
  \draw[->,dashed] (l2box) -- (l3box);
  \begin{scope}[transform canvas={yshift=7pt}]
    \draw[->,dashed] (l2box) -- (l3box);
  \end{scope}
  \begin{scope}[transform canvas={yshift=15pt}]
    \draw[->,dashed] (l2box) -- node[pos=0.5, above] (d2) {$0.5$} (l3box);
  \end{scope}
  \node[above of = d2, node distance = 25pt] {relu};
  \begin{scope}[transform canvas={yshift=-7pt}]
    \draw[->,dashed] (l2box) -- (l3box);
  \end{scope}
  \begin{scope}[transform canvas={yshift=-15pt}]
    \draw[->,dashed] (l2box) -- (l3box);
  \end{scope}
  \node[right of = l3top, node distance = 65pt] (l4top) {};
  \node[right of = l3bottom, node distance = 65pt] (l4bottom) {};
  \node[right of = l3, node distance = 65pt] (l4) {$128 \to 15$};
  \node[above of = l4top, node distance = 15pt] {Dense};
  \node[draw, fit=(l4top)(l4bottom)(l4)] (l4box) {};
  \draw[->,dashed] (l3box) -- (l4box);
  \begin{scope}[transform canvas={yshift=7pt}]
    \draw[->,dashed] (l3box) -- (l4box);
  \end{scope}
  \begin{scope}[transform canvas={yshift=15pt}]
    \draw[->,dashed] (l3box) -- node[pos=0.5, above] (d3) {$0.5$} (l4box);
  \end{scope}
  \node[above of = d3, node distance = 25pt] {relu};
  \begin{scope}[transform canvas={yshift=-7pt}]
    \draw[->,dashed] (l3box) -- (l4box);
  \end{scope}
  \begin{scope}[transform canvas={yshift=-15pt}]
    \draw[->,dashed] (l3box) -- (l4box);
  \end{scope}
  \node[right of=l4, node distance = 65pt] (output) {$\mathbf{x_0}$};
  \draw[->] (l4box) -- (output);
\end{tikzpicture}
  \caption{Initialization network architecture for the inflow forecaster. Activation functions and dropout values between layers are given above the neuron links. The seasonal parameter is given by the current week.}
  \label{fig:flowinitarch}
\end{figure}

The inflow forecaster is slightly more involved. The reason for this is that we wish to forecast inflows to all 15 power stations simultaneously. Besides, inflows between adjacent stations in the river are positively correlated. The time increment is days as opposed to hours, and the forecast horizon is set to one week. The initializer architecture is shown in Figure~\ref{fig:flowinitarch}. It has the same structure as the price network, but uses the current week as the seasonal input parameter.

\begin{figure}
  \centering
  \begin{tikzpicture}
  \node at (0,0) (anchor) {};
  \node[left of=anchor,xshift=-100pt] (day) {$t$};
  \node[above of = day, node distance = 10pt] (state) {$\mathbf{x_{t-1}}$};
  \node[below of = day, node distance = 10pt] (week) {$w$};
  \node[right of = state, node distance = 65pt] (l1top) {};
  \node[right of = week, node distance = 65pt] (l1bottom) {};
  \node[right of = day, node distance = 60pt] (l1) {$17 \to 128$};
  \node[draw,fit=(l1top)(l1bottom)(l1)] (l1box) {};
  \node[above of = l1top, node distance = 15pt] {GRU};
  \draw[->] (state) -- node[pos = 0.5] (sre) {} ($(l1box) + (-27pt,10pt)$);
  \draw[->] (day) -- ($(l1box) + (-27pt,0)$);
  \draw[->] (week) -- ($(l1box) + (-27pt,-10pt)$);
  \node[right of = l1top, node distance = 65pt] (l2top) {};
  \node[right of = l1bottom, node distance = 65pt] (l2bottom) {};
  \node[right of = l1, node distance = 65pt] (l2) {$128 \to 128$};
  \node[draw, fit=(l2top)(l2bottom)(l2)] (l2box) {};
  \node[above of = l2top, node distance = 15pt] {GRU};
  \draw[->,dashed] (l1box) -- (l2box);
  \begin{scope}[transform canvas={yshift=7pt}]
    \draw[->,dashed] (l1box) -- (l2box);
  \end{scope}
  \begin{scope}[transform canvas={yshift=15pt}]
    \draw[->,dashed] (l1box) -- node[pos=0.5, above] (d) {$0.4$} (l2box);
    \node[above of = d, node distance = 10pt] {relu};
  \end{scope}
  \begin{scope}[transform canvas={yshift=-7pt}]
    \draw[->,dashed] (l1box) -- (l2box);
  \end{scope}
  \begin{scope}[transform canvas={yshift=-15pt}]
    \draw[->,dashed] (l1box) -- (l2box);
  \end{scope}
  \node[right of = l2top, node distance = 65pt] (l3top) {};
  \node[right of = l2bottom, node distance = 65pt] (l3bottom) {};
  \node[right of = l2, node distance = 65pt] (l3) {$128 \to 15$};
  \node[draw, fit=(l3top)(l3bottom)(l3)] (l3box) {};
  \draw[->] (l2box) -- (l3box);
  \begin{scope}[transform canvas={yshift=7pt}]
    \draw[->] (l2box) -- (l3box);
  \end{scope}
  \begin{scope}[transform canvas={yshift=15pt}]
    \draw[->] (l2box) -- node[pos=0.5, above] {relu} (l3box);
  \end{scope}
  \begin{scope}[transform canvas={yshift=-7pt}]
    \draw[->] (l2box) -- (l3box);
  \end{scope}
  \begin{scope}[transform canvas={yshift=-15pt}]
    \draw[->] (l2box) -- (l3box);
  \end{scope}
  \node[right of=l3, node distance = 60pt] (output) {$\mathbf{x_t}$};
  \draw[->] (l3box) -- node[pos = 0.5] (sr0) {} (output);
  \node[above of = sr0, node distance = 35pt] (sr1) {};
  \node[above of = sre] (sr2) {};
  \draw (sr0.center) -- (sr1.center) -- (sr2.south) -- (sre.center);
\end{tikzpicture}
  \caption{Sequence generation network for the inflow forecaster. The output is fed back into the system until the full sequence is generated. The seasonal parameters are given by the current week and day.}
  \label{fig:flowseqarch}
\end{figure}
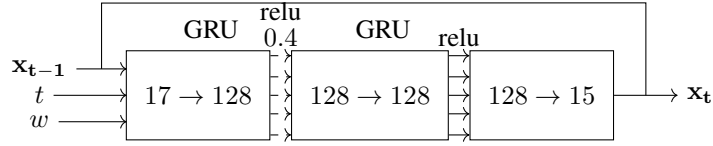

The sequence generator is a three-layered network with larger hidden layers than the price forecaster. The first layer is a GRU based layer with $128$ hidden neurons. It is followed by a dense layer of $128$ hidden neurons and ReLU activation functions. A dropout layer of value $0.4$ is applied between the recurrent layer and the dense layer. The output layer combines the neuron values of the dense layer into the next sequence value, which now has length $15$. Again, the computed flow vector is fed back into the network until inflows for the full week have been computed. Apart from the previous sequence value, the sequence generator is also supplied with the current week and the current day. The sequence generator architecture for the flow forecaster is presented in Fig.~\ref{fig:flowseqarch}.

\subsubsection{Training}
\label{sec:training}

The training procedure involves standard methods from the field of deep learning. The price data is sorted into daily curves, and the inflow data is sorted into weekly curves. Appropriate seasonal parameters are added to each data chunk. About $\SI{10}{\percent}$ of each dataset is removed and used as validation sets. The training procedure operates in epochs. Each epoch, the network generates the desired output using the given input in each chunk and then compares it to the actual output in the chunk using a mean-square cost function. The derivatives of the cost function are then propagated backward through the network using the chain rule and are then used to update the network weights. We use the ADAM~\cite{adam} optimizer to improve the weights, which is known to be effective when training networks on sequential data. After updating all weights, we reset the internal memory of the recurrent layers. This is done to ensure that the performance of the network does not rely on the order it visits the datasets. We also randomize the order of the chunks of each epoch. Every fifth epoch, the performance is evaluated on the validation set, which has not been used to train the weights. If the validation error is increasing, while the training performance still decreases, we stop the training procedure. This is a standard approach to reduce overfitting. Most techniques required for implementing this training procedure are available as tools in Flux.jl.

\subsection{Day-ahead model}
\label{sec:day-ahead-model}

The day-ahead model used in this work is similar to the model introduced in~\cite{Fleten2007}. We repeat the general structure and also highlight the key details of our model. Further, we sketch how the model is implemented in our software framework.

\subsubsection{General structure}
\label{sec:general-structure}

In the first stage, we model the day-ahead orders to be submitted to the market. Let $\mathcal{T} = \braces{t_1,\dots,t_{24}}$ denote indices for the $24$-hour horizon of the upcoming day. The set $\mathcal{B} = \braces{b_1,\dots,b_B}$ make up blocks of consecutive hours in the $24$-hour period. In order to avoid non-linear relations in the model, we fix a set of hourly price levels $\mathcal{P}_t = \braces{p_{1,t},\dots,p_{P,t}}$ to bid at beforehand. We explain how these prices are chosen in a following subsection. The block order prices are determined by calculating averages of the available prices levels over the given blocks. We introduce $x^I_t, x^D_{p,t}$, and $x^B_{p,b}$ to represent price independent orders, price dependent orders, and block orders respectively. As per NordPool regulations, the volumes in a price dependent sell order have to be constant or increasing with increasing prices. We enforce this using constraints. In addition, we constrain the total volume offered to the market to not exceed $200\%$ of the production capacity. Consequently, we allow imbalances in the order commitments, but limit the maximum imbalance already in the first stage.

In the second stage, we model the order commitments after price realization as well as the production schedule after inflow realization. We introduce the random variables $\rho^{\omega}_t$, that describe the hourly market price. Let $y^H_t$ and $y^B_b$ represent the committed hourly volumes and the committed block volumes respectively. Every hour $t$, the dispatched hourly volumes are determined through linear interpolation:
\begin{equation*}
  y^H_t = x^I_t + \frac{\rho^{\omega}_t-p_{i,t}}{p_{i+1,t} - p_{i,t}}x^D_{i+1,t} + \frac{p_{i+1,t}-\rho^{\omega}_t}{p_{i+1,t} - p_{i,t}}x^D_{i, t} \quad p_{i,t} \leq \rho^{\omega}_t \leq p_{i+1,t}.
\end{equation*}
The dispatched block volumes are given by
\begin{equation*}
  y^B_b = \sum_{\mathclap{p : \bar{p}(p,b) \leq \bar{\rho}^{\omega}_b}} x^B_{p,b}
\end{equation*}
where
\begin{equation*}
  \bar{p}(i,b) = \frac{1}{|b|}\sum_{t \in b} p_{i,t}
\end{equation*}
and
\begin{equation*}
  \bar{\rho}^{\omega}_b = \frac{1}{|b|}\sum_{t \in b} \rho^{\omega}_t.
\end{equation*}
Next, we model the production. Let $\mathcal{H} = \braces{h_1,\dots,h_{15}}$ index the 15 hydroelectric power stations in Skellefteälven. For each plant and hour, let $Q_{h,s,t}$ and $S_{h,t}$ denote the water discharged and spilled, respectively. Further, let $P_t$ denote the total volume of electricity produced each hour. We employ a piecewise linear approximation of the generation curve of each station. In other words,
\begin{equation*}
  P_t = \sum_{s \in \mathcal{S}} \mu_{h,s}Q_{h,s,t},
\end{equation*}
where $\mu_{h,s}$ is the marginal production equivalent of station $h$ and segment $s$. The load balance is given by
\begin{equation*}
  y^H + \sum_{t \in b, b \in \mathcal{B}} y^B_b  - P_t = y^{+}_t - y^{-}_t.
\end{equation*}
In each hour, any imbalance between committed volumes and produced volumes is equal to the difference between the imbalance variables $y^{+}_t$ and $y^{-}_t$. Any shortage $y^{+}_t$ is bought from the balancing market, and any surplus $y^{-}_t$ is sold to the balancing market. Finally, let $R_{h,t}$ denote the reservoir contents in plant $h$ during hour $t$. Flow conservation each hour is given by
\begin{equation*}
  \begin{aligned}
    M_{h,t} = \;&M_{h,t-1} \\
    &+ \sum_{i \in \mathcal{Q}_u(h)} \sum_{s \in \mathcal{S}} Q_{i,s,t-\tau_{ih}} + \sum_{i \in \mathcal{S}_u(h)} S_{i,t-\tau_{ih}} + V^{\omega}_h \\
    &- \sum _{s \in \mathcal{S}} Q_{h,s,t} - S_{h,t}
  \end{aligned}
\end{equation*}
Here, $V^{\omega}_h$ are random variables describing the local inflow to each plant. The sets $\mathcal{Q}_u(h)$ and $\mathcal{S}_u(h)$ contain upstream plants where discharge and/or spillage can reach plant $h$ through connecting waterways. Note that the water travel times $\tau_{ih}$ between power stations are included in the incoming flow to each plant. Internally, this is modeled by introducing auxiliary variables and constraints. Variable limits and the introduced parameters are all included in the deterministic data sets for Skellefteälven. The revenue from a production schedule satisfying the above relations is given by
\begin{equation*}
  \sum_{t \in \mathcal{T}} \rho^{\omega}_ty^H_t + \sum_{b \in \mathcal{B}}\abs{b}\bar{\rho}^{\omega}_by^B_b + \sum_{t \in \mathcal{T}}\alpha_t\rho^{\omega}_ty^{-}_t - \beta_t\rho^{\omega}_ty^{+}_t + W(M_{1,24},\dots,M_{15,24}).
\end{equation*}
Note that, for any committed block order $y^B_b$, the order volume is dispatched every hour in the block at average market price. Hence, $\abs{b}\bar{\rho}^{\omega}_by^B_b$ is earned. The imbalance volumes are traded at penalized prices, using penalty factors $\alpha_t$ and $\beta_t$, for discouragement. It is hard to accurately model this penalty. Here, we use a $15\%$ penalty during peak hours, and $10\%$ otherwise. These values are based on observations of historic values, but can not be considered accurate. The final term in the revenue is the expected future value of water, which is a function of the water volumes that remain in the reservoirs after the period. In the following section, we introduce a polyhedral approximation of this function that can be modelled with linear terms. For now, we simply denote the water value by $W$. In summary, stochastic program modeling the day-ahead problem is in essence given by
\begin{equation} \label{eq:dayahead}
  \begin{aligned}
   \maximize_{x^I_t, x^D_{i,h}, x^B_{i,b}} & \quad \expect[\omega]{\sum_{t \in \mathcal{T}} \rho^{\omega}_ty^H_t + \sum_{b \in \mathcal{B}}\abs{b}\bar{\rho}^{\omega}_by^B_b + \sum_{t \in \mathcal{T}}\alpha_t\rho^{\omega}_ty^{-}_t - \beta_t\rho^{\omega}_ty^{+}_t + W} \\
   \st & \quad y^H_t = x^I_t + \frac{\rho^{\omega}_t-p_{i,t}}{p_{i+1,t} - p_{i,t}}x^D_{i+1,t} + \frac{p_{i+1,t}-\rho^{\omega}_t}{p_{i+1,t} - p_{i,t}}x^D_{i, t} \\
   & \quad y^B_b = \sum_{p : \bar{p}(p,b) \leq \bar{\rho}^{\omega}_b} x^B_{p,b} \\
   & \quad P_t = \sum_{s \in \mathcal{S}} \mu_{h,s}Q_{h,s,t} \\
   & \quad y^H + \sum_{\mathclap{t \in b, b \in \mathcal{B}}} y^B_b  - P_t = y^{+}_t - y^{-}_t \\
   & \quad \begin{aligned}
     M_{h,t} = \;&M_{h,t-1} \\
     &+ \sum_{i \in \mathcal{Q}_u(h)} \sum_{s \in \mathcal{S}} Q_{i,s,t-\tau_{ih}} + \sum_{i \in \mathcal{S}_u(h)} S_{i,t-\tau_{ih}} + V^{\omega}_h \\
     &- \sum _{s \in \mathcal{S}} Q_{h,s,t} - S_{h,t} \\
   \end{aligned} \\
   & \quad 0 \leq Q_{h,s,t} \leq \bar{Q}_{h,s} \\
   & \quad 0 \leq M_{h,t} \leq \bar{M}_h \\
   & \quad y^H_t \geq 0,\; y^B_b \geq 0,\; y^{+}_t \geq 0,\; y^{-}_t \geq 0 \\
   & \quad P_t \geq 0,\; S_{h,t} \geq 0.
 \end{aligned}
\end{equation}

\subsubsection{Water evaluation}
\label{sec:water-evaluation}

The expected value of saving water has a large impact on the second-stage production schedule. If the water value is large, then it could be optimal to not produce, settle committed orders in the balancing market, and save water. Likewise, if the water value is small, it could be optimal to overproduce and sell the excess in the balancing market. Consequently, the water value will evidently also impact the optimal order strategy because the order commitments are instrumental in both scenarios. Thus, the accuracy of the water evaluation is critical for hydropower producers participating in the day-ahead market. A naive approach is to assume that all excess water can be used to produce electricity sold at some expected future price. This leads to coarse order strategies as the optimal strategy is governed by price variations around the expected future price. We adopt a slightly more involved approach.

Consider a dummy stochastic program, where the first-stage decision is to decide the reservoir contents of every power station before the upcoming week. Next, we realize a weekly sequence of inflows and daily price curves, using the noise-driven RNN forecasters, and optimize the weekly production of energy sold to the market. This simplified week-ahead problem is given by

\begin{equation} \label{eq:weekahead}
  \begin{aligned}
   \maximize_{M_{h,0}} & \quad \expect[\omega]{\sum_{t \in \tilde{\mathcal{T}}} \rho^{\omega}_tP_t} \\
   \st & \quad P_t = \sum_{s \in \mathcal{S}} \mu_{h,s}Q_{h,s,t} \\
   & \quad \begin{aligned}
    M_{h,t} &= M_{h,t-1} \\
            &+ \sum_{i \in \mathcal{Q}_u(h)} \sum_{s \in \mathcal{S}} Q_{i,s,t-\tau_{ih}} + \sum_{i \in \mathcal{S}_u(h)} S_{i,t-\tau_{ih}} + V^{\omega}_h \\
            &- \sum _{s \in \mathcal{S}} Q_{h,s,t} - S_{h,t} \\
          \end{aligned} \\
   & \quad 0 \leq Q_{h,s,t} \leq \bar{Q}_{h,s} \\
   & \quad 0 \leq M_{h,t} \leq \bar{M}_h \\
   & \quad P_t \geq 0,\; S_{h,t} \geq 0,
 \end{aligned}
\end{equation}
\noindent
where the time-horizon $\tilde{\mathcal{T}} = \braces{t_1,\dots,t_{168}}$ is now a week. The objective function $W = \expect[\omega]{W^{\omega}(M_{1,0},\dots,M_{15,0})}$ of this problem will be used as a water value function. The problem~\eqref{eq:weekahead} is trivial since the optimal decision will be to fill the reservoirs with enough water to be able to run at maximum capacity in the worst-case scenario. However, information about the water value can be extracted by solving~\eqref{eq:weekahead} with an L-shaped type method. The L-shaped method generates cutting planes of the form
\begin{equation} \label{eq:wcuts}
  \sum_{h \in \mathcal{H}}\partial W_{c,h}M_{h,0}+ W \geq w_c.
\end{equation}
\noindent
This form supports for the concave objective function $W$, which is a function of reservoir content in the system. Hence, after the algorithm has converged we have access to a polyhedral approximation of $W$ in the form of a collection of such cuts as~\eqref{eq:wcuts}. We can use these cuts to put an approximate value of the remaining volumes of water present in the reservoirs after meeting order commitments. The water value approximation enters the day-ahead problem~\eqref{eq:dayahead} in the following way:
\begin{equation*}
\begin{aligned}
 \maximize_{} & \quad \dots + W \\
 \st & \quad \vdots \\
 & \quad \sum_{h \in \mathcal{H}}\partial W_{c,h}M_{h,24}+ W \geq w_c \quad c \in \mathcal{C}. \\
\end{aligned}
\end{equation*}
In practice, we use a multiple-cut formulation
\begin{equation*}
  W = \sum_{i = 1}^{N} W_i
\end{equation*}
as the L-shaped method solves the week-ahead problem with a large number of scenarios $N$ more efficiently in this way. The end result is still a collection of cuts that approximate a polyhedral water value function of the final reservoir volumes.

\subsubsection{Price levels}
\label{sec:pricelevels}

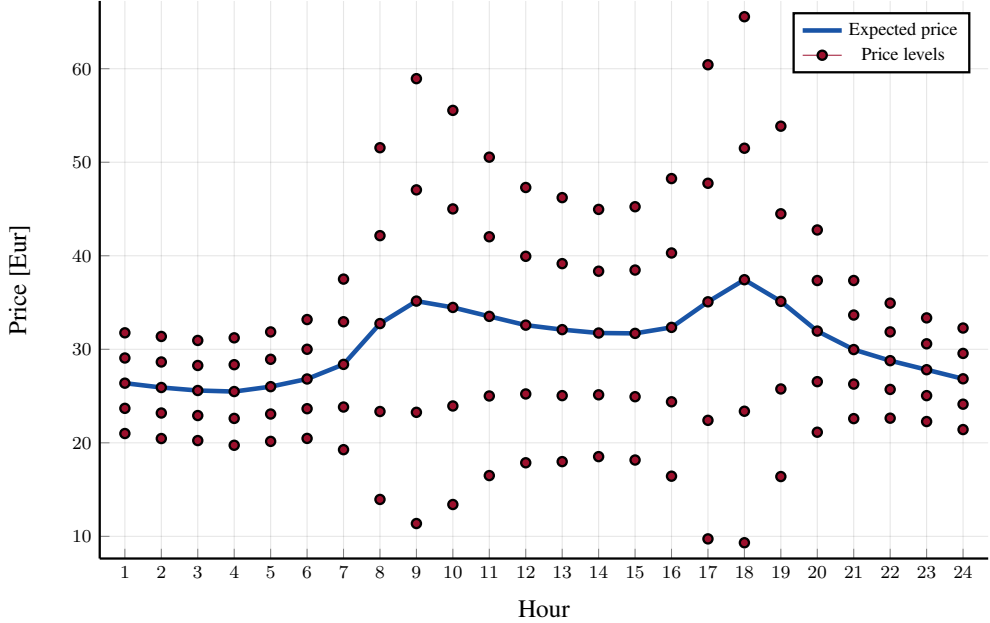
\begin{figure}
  \centering
  \resizebox{0.8\textwidth}{!}{
\begin{tikzpicture}[]
\begin{axis}[height = {101.6mm}, ylabel = {Price [Eur]}, xmin = {0.31000000000000005}, xmax = {24.69}, ymax = {67.24766457634726}, xlabel = {Hour}, unbounded coords=jump,scaled x ticks = false,xlabel style = {font = {\fontsize{11 pt}{14.3 pt}\selectfont}, color = {rgb,1:red,0.00000000;green,0.00000000;blue,0.00000000}, draw opacity = 1.0, rotate = 0.0},xmajorgrids = true,xtick = {1.0,2.0,3.0,4.0,5.0,6.0,7.0,8.0,9.0,10.0,11.0,12.0,13.0,14.0,15.0,16.0,17.0,18.0,19.0,20.0,21.0,22.0,23.0,24.0},xticklabels = {$1$,$2$,$3$,$4$,$5$,$6$,$7$,$8$,$9$,$10$,$11$,$12$,$13$,$14$,$15$,$16$,$17$,$18$,$19$,$20$,$21$,$22$,$23$,$24$},xtick align = inside,xticklabel style = {font = {\fontsize{8 pt}{10.4 pt}\selectfont}, color = {rgb,1:red,0.00000000;green,0.00000000;blue,0.00000000}, draw opacity = 1.0, rotate = 0.0},x grid style = {color = {rgb,1:red,0.00000000;green,0.00000000;blue,0.00000000},
draw opacity = 0.1,
line width = 0.5,
solid},axis x line* = left,x axis line style = {color = {rgb,1:red,0.00000000;green,0.00000000;blue,0.00000000},
draw opacity = 1.0,
line width = 1,
solid},scaled y ticks = false,ylabel style = {font = {\fontsize{11 pt}{14.3 pt}\selectfont}, color = {rgb,1:red,0.00000000;green,0.00000000;blue,0.00000000}, draw opacity = 1.0, rotate = 0.0},ymajorgrids = true,ytick = {10.0,20.0,30.0,40.0,50.0,60.0},yticklabels = {$10$,$20$,$30$,$40$,$50$,$60$},ytick align = inside,yticklabel style = {font = {\fontsize{8 pt}{10.4 pt}\selectfont}, color = {rgb,1:red,0.00000000;green,0.00000000;blue,0.00000000}, draw opacity = 1.0, rotate = 0.0},y grid style = {color = {rgb,1:red,0.00000000;green,0.00000000;blue,0.00000000},
draw opacity = 0.1,
line width = 0.5,
solid},axis y line* = left,y axis line style = {color = {rgb,1:red,0.00000000;green,0.00000000;blue,0.00000000},
draw opacity = 1.0,
line width = 1,
solid},    xshift = 0.0mm,
    yshift = 0.0mm,
    axis background/.style={fill={rgb,1:red,1.00000000;green,1.00000000;blue,1.00000000}}
,legend style = {color = {rgb,1:red,0.00000000;green,0.00000000;blue,0.00000000},
draw opacity = 1.0,
line width = 1,
solid,fill = {rgb,1:red,1.00000000;green,1.00000000;blue,1.00000000},font = {\fontsize{8 pt}{10.4 pt}\selectfont}},colorbar style={title=}, ymin = {7.628412739936932}, width = {152.4mm}]\addplot+ [color = {rgb,1:red,0.09803922;green,0.32941176;blue,0.65098039},
draw opacity = 1.0,
line width = 2,
solid,mark = none,
mark size = 2.0,
mark options = {
    color = {rgb,1:red,0.00000000;green,0.00000000;blue,0.00000000}, draw opacity = 1.0,
    fill = {rgb,1:red,0.09803922;green,0.32941176;blue,0.65098039}, fill opacity = 1.0,
    line width = 1,
    rotate = 0,
    solid
}]coordinates {
(1.0, 26.377586431503296)
(2.0, 25.91467671966553)
(3.0, 25.592450896263124)
(4.0, 25.481782537460326)
(5.0, 26.004895790100097)
(6.0, 26.825305461883545)
(7.0, 28.385489170074464)
(8.0, 32.749079216003416)
(9.0, 35.15360792827606)
(10.0, 34.47467072868347)
(11.0, 33.52035451698303)
(12.0, 32.583070892333986)
(13.0, 32.102346574783326)
(14.0, 31.740470558166503)
(15.0, 31.701566961288453)
(16.0, 32.342629470825194)
(17.0, 35.07680814361572)
(18.0, 37.43803865814209)
(19.0, 35.123736673355104)
(20.0, 31.9463049659729)
(21.0, 29.970932359695436)
(22.0, 28.784440793991088)
(23.0, 27.81591905593872)
(24.0, 26.84460446166992)
};
\addlegendentry{Expected price}
\addplot+[draw=none, color = {rgb,1:red,0.61568627;green,0.06274510;blue,0.17647059},
draw opacity = 1.0,
line width = 0,
solid,mark = *,
mark size = 2.0,
mark options = {
    color = {rgb,1:red,0.00000000;green,0.00000000;blue,0.00000000}, draw opacity = 1.0,
    fill = {rgb,1:red,0.61568627;green,0.06274510;blue,0.17647059}, fill opacity = 1.0,
    line width = 1,
    rotate = 0,
    solid
}] coordinates {
(1.0, 20.99341396253017)
(2.0, 20.45733220765932)
(3.0, 20.239489859538217)
(4.0, 19.737804863058187)
(5.0, 20.152802829257215)
(6.0, 20.470774425740707)
(7.0, 19.26656043356128)
(8.0, 13.945209308787053)
(9.0, 11.369126198144443)
(10.0, 13.403605200008595)
(11.0, 16.49702675391818)
(12.0, 17.865773269678996)
(13.0, 17.988791440942336)
(14.0, 18.523584451985485)
(15.0, 18.15928997909782)
(16.0, 16.431537024584017)
(17.0, 9.728789733695685)
(18.0, 9.315750056061752)
(19.0, 16.391890863433396)
(20.0, 21.133446449558875)
(21.0, 22.586666222843522)
(22.0, 22.63350740400468)
(23.0, 22.273971688418406)
(24.0, 21.415767819469757)
};
\addlegendentry{Price levels}
\addplot+[draw=none, color = {rgb,1:red,0.61568627;green,0.06274510;blue,0.17647059},
draw opacity = 1.0,
line width = 0,
solid,mark = *,
mark size = 2.0,
mark options = {
    color = {rgb,1:red,0.00000000;green,0.00000000;blue,0.00000000}, draw opacity = 1.0,
    fill = {rgb,1:red,0.61568627;green,0.06274510;blue,0.17647059}, fill opacity = 1.0,
    line width = 1,
    rotate = 0,
    solid
},forget plot] coordinates {
(1.0, 23.685500197016733)
(2.0, 23.186004463662425)
(3.0, 22.91597037790067)
(4.0, 22.609793700259257)
(5.0, 23.078849309678656)
(6.0, 23.648039943812126)
(7.0, 23.826024801817873)
(8.0, 23.347144262395233)
(9.0, 23.26136706321025)
(10.0, 23.93913796434603)
(11.0, 25.008690635450606)
(12.0, 25.22442208100649)
(13.0, 25.04556900786283)
(14.0, 25.132027505075996)
(15.0, 24.930428470193135)
(16.0, 24.387083247704606)
(17.0, 22.402798938655703)
(18.0, 23.37689435710192)
(19.0, 25.75781376839425)
(20.0, 26.539875707765887)
(21.0, 26.27879929126948)
(22.0, 25.708974098997885)
(23.0, 25.044945372178564)
(24.0, 24.130186140569837)
};
\addplot+[draw=none, color = {rgb,1:red,0.61568627;green,0.06274510;blue,0.17647059},
draw opacity = 1.0,
line width = 0,
solid,mark = *,
mark size = 2.0,
mark options = {
    color = {rgb,1:red,0.00000000;green,0.00000000;blue,0.00000000}, draw opacity = 1.0,
    fill = {rgb,1:red,0.61568627;green,0.06274510;blue,0.17647059}, fill opacity = 1.0,
    line width = 1,
    rotate = 0,
    solid
},forget plot] coordinates {
(1.0, 26.377586431503296)
(2.0, 25.91467671966553)
(3.0, 25.592450896263124)
(4.0, 25.481782537460326)
(5.0, 26.004895790100097)
(6.0, 26.825305461883545)
(7.0, 28.385489170074464)
(8.0, 32.749079216003416)
(9.0, 35.15360792827606)
(10.0, 34.47467072868347)
(11.0, 33.52035451698303)
(12.0, 32.583070892333986)
(13.0, 32.102346574783326)
(14.0, 31.740470558166503)
(15.0, 31.701566961288453)
(16.0, 32.342629470825194)
(17.0, 35.07680814361572)
(18.0, 37.43803865814209)
(19.0, 35.123736673355104)
(20.0, 31.9463049659729)
(21.0, 29.970932359695436)
(22.0, 28.784440793991088)
(23.0, 27.81591905593872)
(24.0, 26.84460446166992)
};
\addplot+[draw=none, color = {rgb,1:red,0.61568627;green,0.06274510;blue,0.17647059},
draw opacity = 1.0,
line width = 0,
solid,mark = *,
mark size = 2.0,
mark options = {
    color = {rgb,1:red,0.00000000;green,0.00000000;blue,0.00000000}, draw opacity = 1.0,
    fill = {rgb,1:red,0.61568627;green,0.06274510;blue,0.17647059}, fill opacity = 1.0,
    line width = 1,
    rotate = 0,
    solid
},forget plot] coordinates {
(1.0, 29.06967266598986)
(2.0, 28.643348975668633)
(3.0, 28.268931414625577)
(4.0, 28.353771374661395)
(5.0, 28.930942270521538)
(6.0, 30.002570979954964)
(7.0, 32.944953538331056)
(8.0, 42.1510141696116)
(9.0, 47.04584879334187)
(10.0, 45.01020349302091)
(11.0, 42.03201839851546)
(12.0, 39.94171970366148)
(13.0, 39.15912414170382)
(14.0, 38.34891361125701)
(15.0, 38.47270545238377)
(16.0, 40.29817569394578)
(17.0, 47.75081734857574)
(18.0, 51.499182959182264)
(19.0, 44.48965957831596)
(20.0, 37.35273422417991)
(21.0, 33.66306542812139)
(22.0, 31.85990748898429)
(23.0, 30.58689273969888)
(24.0, 29.559022782770004)
};
\addplot+[draw=none, color = {rgb,1:red,0.61568627;green,0.06274510;blue,0.17647059},
draw opacity = 1.0,
line width = 0,
solid,mark = *,
mark size = 2.0,
mark options = {
    color = {rgb,1:red,0.00000000;green,0.00000000;blue,0.00000000}, draw opacity = 1.0,
    fill = {rgb,1:red,0.61568627;green,0.06274510;blue,0.17647059}, fill opacity = 1.0,
    line width = 1,
    rotate = 0,
    solid
},forget plot] coordinates {
(1.0, 31.761758900476423)
(2.0, 31.372021231671738)
(3.0, 30.94541193298803)
(4.0, 31.225760211862465)
(5.0, 31.85698875094298)
(6.0, 33.17983649802638)
(7.0, 37.50441790658765)
(8.0, 51.552949123219776)
(9.0, 58.93808965840768)
(10.0, 55.54573625735834)
(11.0, 50.54368228004789)
(12.0, 47.300368514988975)
(13.0, 46.215901708624315)
(14.0, 44.95735666434752)
(15.0, 45.243843943479085)
(16.0, 48.25372191706637)
(17.0, 60.42482655353576)
(18.0, 65.56032726022244)
(19.0, 53.85558248327681)
(20.0, 42.759163482386924)
(21.0, 37.35519849654735)
(22.0, 34.9353741839775)
(23.0, 33.357866423459036)
(24.0, 32.273441103870084)
};
\end{axis}
\end{tikzpicture}}
  \caption{Expected daily electricity price out of $1000$ samples from the RNN forecaster. Two standard deviations above and below the expected price is shown each hour.}
  \label{fig:pricelevels}
\end{figure}

The price-dependent hourly orders and the block orders are specified at pre-chosen price levels. For flexibility, we allow these levels to vary with time. The set of price levels $\mathcal{P}_t$ for each hour is determined using the price forecaster. We sample a large number of scenarios and compute the hourly mean price and the hourly standard deviation of the resulting price curves. Next, in each hour $t$, we define four price levels around the mean price using multiples of the standard deviation for that hour over all the sampled curves. We also include the mean price as an available price level. A set of hourly price levels generated using this method is shown in Fig.~\ref{fig:pricelevels}. The aim of this price level generation method is to allow for flexibility in order placement during hours of large historical price variation. For each block $b \in \mathcal{B}$, we define five possible block order prices $\braces{p_{i,b}}_{i = 1}^{5}$ by computing mean price levels over the hours $t \in b$ for each of the five pre-computed price levels to reflect the market settlement of block orders.

\subsubsection{Model definition}
\label{sec:model-definition}

\sloppy
The day-ahead model is formulated in \jlinl{StochasticPrograms.jl}~\cite{spjl} (SPjl), our general purpose software framework for stochastic programming implemented in the Julia programming language. Optimization models in SPjl are defined using the algebraic modeling language JuMP~\cite{Dunning2017}. To increase readability, we present an abridged version of the day-ahead model implementation in SPjl, by obfuscating parts of the code and making slight syntax changes. The full unabridged model is available at Github~\footnote{\url{https://github.com/martinbiel/HydroModels.jl}}.

\sloppy
First, we define a data structure to describe the uncertain parameters using the \jlinl{@scenario} command. We also create a sampler object, using the \jlinl{@sampler} command, which utilizes the noise-driven RNN forecasters to generate price curves and inflows. The code is shown in Listing~\ref{lst:scenariodef}. Because we want to make use of the forecasters' seasonal capabilities, we also include a date field in the sampler object. The forecasters use the provided date to determine seasonal parameter inputs to the neural networks.

\begin{lstlisting}[language = julia, float, caption = {Day-ahead scenario definition in SPjl}, label = {lst:scenariodef}]
@scenario DayAheadScenario = begin
    ρ::PriceCurve{Float64}
    V::Inflows{Float64}
end

@sampler RNNDayAheadSampler = begin
    date::Date
    price_forecaster::PriceForecaster
    flow_forecaster::FlowForecaster

    @sample DayAheadScenario begin
        price_curve = forecast(sampler.price_forecaster, month(sampler.date))
        flows = forecast(sampler.flow_forecaster, week(sampler.date))
        return DayAheadScenario(PriceCurve(price_curve), Inflows(flows))
    end
end
\end{lstlisting}
\noindent
\sloppy
The day-ahead model definition in SPjl is presented in Listing~\ref{lst:dayaheaddef}. Internally, the \jlinl{@stochastic_model} block creates model recipes for the stage problems without actually instantiating any optimization problems. There are two special lines in Listing~\ref{lst:dayaheaddef}. The \jlinl{@decision} annotation tells the system how the first and second stages are linked. Similarly, the \jlinl{@uncertain} annotation describes how scenario data enters the second stage. A specific second-stage instance is only created when a \jlinl{DayAheadScenario} is supplied. This deferred model instantiation technique allows us to efficiently instantiate subproblems on remote compute nodes, by only passing model recipes and sampler objects with low memory footprint. Consequently, we can solve large-scale day-ahead instances using parallel decomposition schemes.

\begin{lstlisting}[language = julia, float, caption = {Day-ahead problem definition in SPjl. The code has been condensed for readability.}, label = {lst:dayaheaddef}]
@stochastic_model begin
    @stage 1 begin
        @parameters horizon indices data
        @unpack hours, plants, bids, blockbids, blocks = indices
        @unpack hydrodata, regulations = data
        @decision(model, xᴵ[t = hours] >= 0)
        @decision(model, xᴰ[p = bids, t = hours] >= 0)
        @decision(model, 0 <= xᴮ[p = blockbids, b = blocks] <= blocklimit)
        # Increasing bid curve
        @constraint(model, bidcurve[p = bids[1:end-1], t = hours],
            xᴰ[p,t] <= xᴰ[p+1,t])
        # Maximal bids ...
    end
    @stage 2 begin
        @parameters horizon indices data
        @unpack hours, plants, segments, blocks = indices
        @unpack hydrodata, water_value, regulations, bidlevels = data
        @uncertain ρ, V from ξ::DayAheadScenario
        @variable(model, yᴴ[t = hours] >= 0)
        @variable(model, yᴮ[b = blocks] >= 0)
        @recourse(model, y⁺[t = hours] >= 0)
        @recourse(model, y⁻[t = hours] >= 0)
        @variable(model, 0 <= Q[h = plants,s = segments,t = hours] <= Q_max[s])
        @variable(model, S[h = plants,t = hours] >= 0)
        @variable(model, 0 <= M[h = plants,t = hours] <= M_max[p])
        @variable(model, W[i = 1:nindices(water_value)])
        @variable(model, H[t = hours] >= 0)
        @expression(model, net_profit,
            sum(ρ[t]*yᴴ[t] for t in hours)
            + sum(|b|*(mean(ρ[hours_per_block[b]])*yᴮ[b] for b in blocks))
        @expression(model, intraday,
            sum(penalty(ξ,t)*y⁺[t] - reward(ξ,t)*y⁻[t] for t in hours))
        @expression(model, value_of_stored_water,
            -sum(W[i] for i in 1:nindices(water_value)))
        @objective(model, Max, net_profit - intraday + value_of_stored_water)
        # Bid-dispatch links
        @constraint(model, hourlybids[t = hours],
            yᴴ[t] == interpolate(ρ[t], bidlevels, xᴰ[t]) + xᴵ[t])
        @constraint(model, bidblocks[b = blocks],
            yᴮ[b] == sum(xᴮ[j,b] for j in accepted_blockorders(b)))
        # Hydrological balance
        @constraint(model, hydro_constraints[h = plants, t = hours],
            M[h,t] == (t > 1 ? M[h,t-1] : M₀[h])
            + sum(Q[i,s,t-τ] for i in intersect(Qu[h], plants), s in segments)
            + sum(S[i,s,t-τ] for i in intersect(Su[h], plants), s in segments)
            + V[h] - sum(Q[h,s,t] for s in segments) - S[h,t])
        # Production
        @constraint(model, production[t = hours],
            H[t] == sum(μ[s]*Q[h,s,t] for h in plants, s in segments))
        # Load balance
        @constraint(model, loadbalance[t = hours],
            yᴴ[t] + sum(yᴮ[b] for b in active(t)) - H[t] == y⁺[t] - y⁻[t])
        # Water travel time ...
        # Polyhedral water value
        @constraint(model, water_value_approximation[c = 1:ncuts(water_value)],
            sum(∂W[c,h]*M[h,T] for h in plants)
            + sum(W[i] for i in cut_indices(c)) >= w[c])
    end
end
\end{lstlisting}

\subsection{Algorithm details}
\label{sec:algorithm-details}

We discuss the algorithmic procedure used to solve the day-ahead problem~\eqref{eq:dayahead}. We will base the discussion around solving a general stochastic program of the form
\begin{equation} \label{eq:generalsp}
  \minimize_{x \in \mathcal{X}} \quad \expect[\omega]{F(x,\xi(\omega))}.
\end{equation}
We use a sampled average approximation (SAA) scheme to solve~\eqref{eq:generaldayahead}. This is a well-established approach (See for example~\cite{saa} for a theoretical introduction and~\cite{saacomp} for an extensive computational study). In short, the algorithm operates by calculating provable confidence intervals around the optimal value (VRP), proceeding until the length of the confidence interval is within some relative tolerance. We then generate the deterministic solution and calculate a confidence interval around the expected value of the expected value solution (EEV). Afterwards, we calculate a confidence interval around the VSS as the difference between the VRP and the EEV.

The confidence interval is calculated as follows. Let $\braces{\mathbf{\xi}^i}_{i = 1}^{M}$ where $\mathbf{\xi}^i = \braces{\xi_j^i}_{j = 1}^{n}$ be $M$ sets of $n$ i.i.d samples from the distribution of $\xi$. For each $i$, consider the sampled problem
\begin{equation} \label{eq:saa}
  z^*_n = \min_{x \in \mathcal{X}} \quad \frac{1}{n}\sum_{j = 1}^{n}F(x,\xi_j^i).
\end{equation}
The authors of~\cite{saa} show that $\expect[\omega]{z^*_n} \leq z^*$ and that $\expect[\omega]{z^*_{n+1}} \geq \expect[\omega]{z^*_n}$. We can use these results to construct lower bounds $L_{n,M}$ that increases with $n$. The $M$ batches of sample sets can be used to construct an unbiased estimator
\begin{equation*}
  \frac{1}{M}\sum_{i = 1}^{M} \min_{x \in \mathcal{X}} \frac{1}{n} \sum_{j = 1}^{n} F(x, \xi_j^i)
\end{equation*}
of $z^*_n$, which is then used as a lower bound estimate. Next, a confidence interval is constructed around the estimate using the sample variance
\begin{equation*}
  \sigma_{M}^2 = \frac{1}{M-1}\sum_{i = 1}^{M} \parentheses*{\min_{x \in \mathcal{X}} \frac{1}{n} \sum_{j = 1}^{n} F(x, \xi_j^i) - L_{n,M}}^2.
\end{equation*}
For practical purposes, $M$ is relatively small, i.e., on the order of $10^1$. Therefore, we construct confidence intervals around $L_{n,M}$ using the $\alpha$-critical value of the $t$-distribution with $M-1$ degrees of freedom, as opposed to using a normal distribution. The approximate $(1-\alpha)$ confidence interval is then given by
\begin{equation*}
\brackets*{L_{n,M} - \frac{t_{\alpha/2,M-1}\sigma_{M}^2}{\sqrt{M}}, L_{n,M} + \frac{t_{\alpha/2,M-1}\sigma_{M}^2}{\sqrt{M}}}.
\end{equation*}
The upper bound is computed from a suboptimal candidate solution $\hat{x}$. It holds that
\begin{equation*}
  z^*(\hat{x}) = \expect[\omega]{F(\hat{x}, \xi(\omega))} \geq z^*
\end{equation*}
for all $\hat{x}$. For a given $\hat{x}$, we estimate $z^*(\hat{x})$ by sampling $T$ batches $\braces{\mathbf{\xi}^i}_{i = 1}^{T}$ where $\mathbf{\xi}^i = \braces{\xi_j^i}_{j = 1}^{N}$ and construct the unbiased estimator
\begin{equation} \label{eq:upperbound}
  \frac{1}{T}\sum_{i = 1}^{M} \frac{1}{N} \sum_{j = 1}^{N} F(\hat{x}, \xi_j^i),
\end{equation}
which is then used as an upper bound $U_{n,T}$ estimate. Here, $N$ can be a large number because~\eqref{eq:upperbound} simply amounts to solving $N$ independent smaller optimization problems. Again, the sample variance and the $t$-distribution are used to calculate a confidence interval around $U_{n,T}$:
\begin{equation*}
  \brackets*{U_{n,T} - \frac{t_{\alpha/2,T-1}\sigma_{T}^2}{\sqrt{T}}, U_{n,T} + \frac{t_{\alpha/2,T-1}\sigma_{T}^2}{\sqrt{T}}}.
\end{equation*}
For convenience, $\hat{x}$ is calculated by solving a single SAA instance of size $n$ at each iteration. Now, the algorithm proceeds by iteratively increasing $n$, calculating a $1-2\alpha$ confidence interval around $U_{n,T}-L_{n,M}$:
\begin{equation*}
  \brackets*{0, U_{n,T} - L_{n,M} + \frac{t_{\alpha/2,T-1}\sigma_{T}^2}{\sqrt{T}} + \frac{t_{\alpha/2,M-1}\sigma_{M}^2}{\sqrt{M}}}
\end{equation*}
and repeating until the gap is within some relative tolerance. It could hold that $U_{n,T} < L_{n,M}$ due to the sampling error, whereupon the procedure simply continues with a slightly larger $n$.

The SAA procedure involves solving numerous sampled instances of increasing size. This is computationally demanding, so we employ parallelization strategies. We distribute every instance on a $\num{32}$-core compute node, using the capabilities of SPjl. Further, we solve the instances using a distributed L-shaped algorithm. A review of our computational experience with these algorithms is given in~\cite{distlshaped} and~\cite{spjl}. Our results indicate that trust-region regularization is appropriate in reducing convergence times when solving day-ahead problem instances. Furthermore, we found that an aggregation scheme can significantly reduce computational times in the distributed setting, through load balancing between master and workers as well as reduced communication latency.

Finally, we sketch the calculation of the value of the stochastic solution. After convergence, the result of the above procedure is a confidence interval $\brackets{L_{VRP},U_{VRP}}$ around the optimal value (VRP) of the day-ahead problem. Next, the expected value solution is cheaply computed through
\begin{equation*}
  \bar{x} = \argmin[x \in \mathcal{X}]{F(x, \bar{\xi})}
\end{equation*}
where $\bar{\xi} = \expect[\omega]{\xi(\omega)}$ is the expected scenario. A confidence interval around the expected result of using the expected value solution (EEV)
\begin{equation*}
  EEV = \frac{1}{\bar{N}}\sum_{i = 1}^{\bar{N}} F(\bar{x}, \xi_i),
\end{equation*}
is given by
\begin{equation*}
  \brackets*{EEV-\frac{z_{\alpha/2}\sigma^2_{EEV}}{\sqrt{\bar{N}}}, EEV+\frac{z_{\alpha/2}\sigma^2_{EEV}}{\sqrt{\bar{N}}}},
\end{equation*}
where
\begin{equation*}
  \sigma_{EEV}^2 = \frac{1}{\bar{N}-1}\sum_{i = 1}^{\bar{N}} \parentheses{F(\bar{x}, \xi_i) - EEV}^2.
\end{equation*}
Here, $\braces{\xi_i}_{i = 1}^{\bar{N}}$ is a large number of sampled scenarios. Now, if there is no overlap between the confidence interval around VRP and the confidence interval around EEV, there is a VSS that is statistically significant to the chosen significance level. A confidence interval around this VSS is then given by
\begin{equation*}
  \brackets*{L_{VRP} - U_{EEV}, U_{VRP} - L_{EEV}}
\end{equation*}

\section{Numerical Experiments}
\label{sec:numerical_experiments}

In this section, we present our experimental results.

\subsection{Day-ahead forecasting}
\label{sec:day-ahead-forec}

\begin{figure}
  \centering
  \input{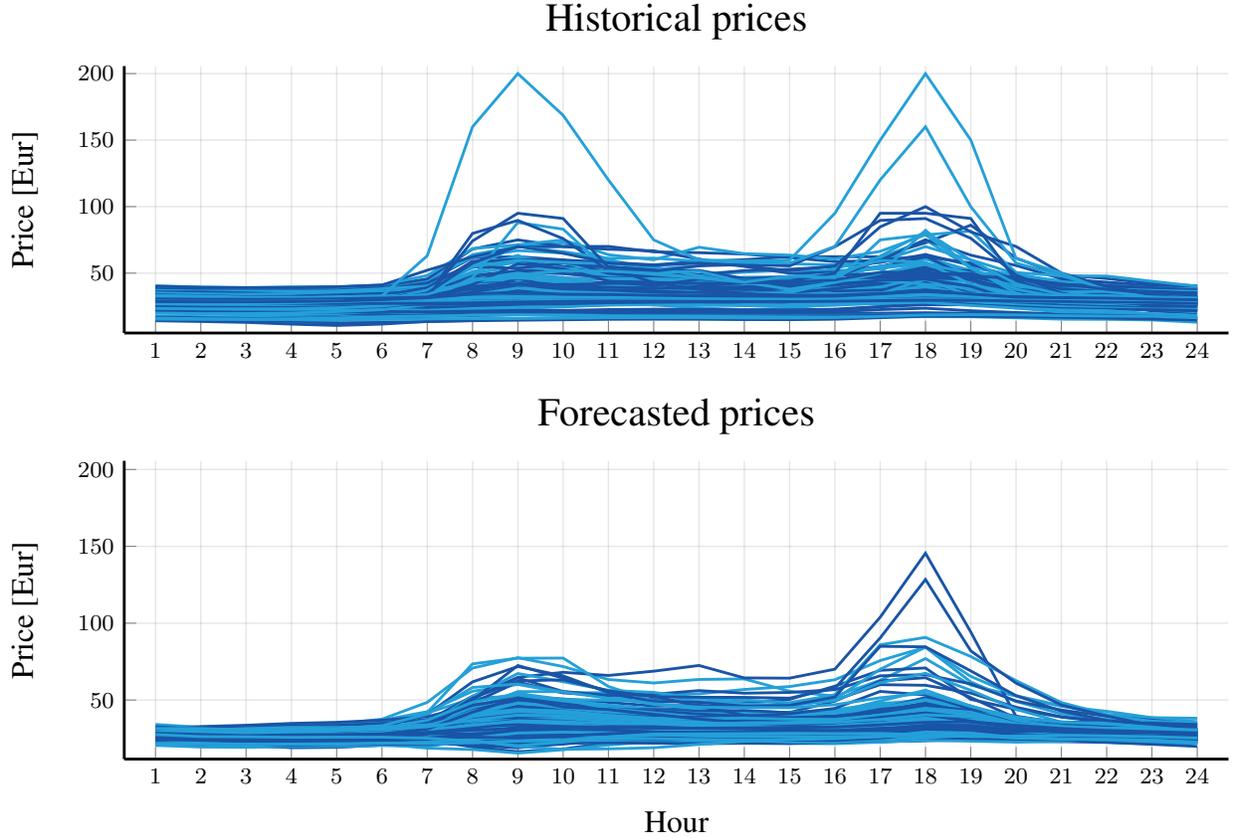}
  \caption{Historical electricity price curves in January and electricity price curves generated using the RNN forecaster in the same period.}
  \label{fig:dayaheadforecast}
\end{figure}

\begin{figure}
  \centering
  \input{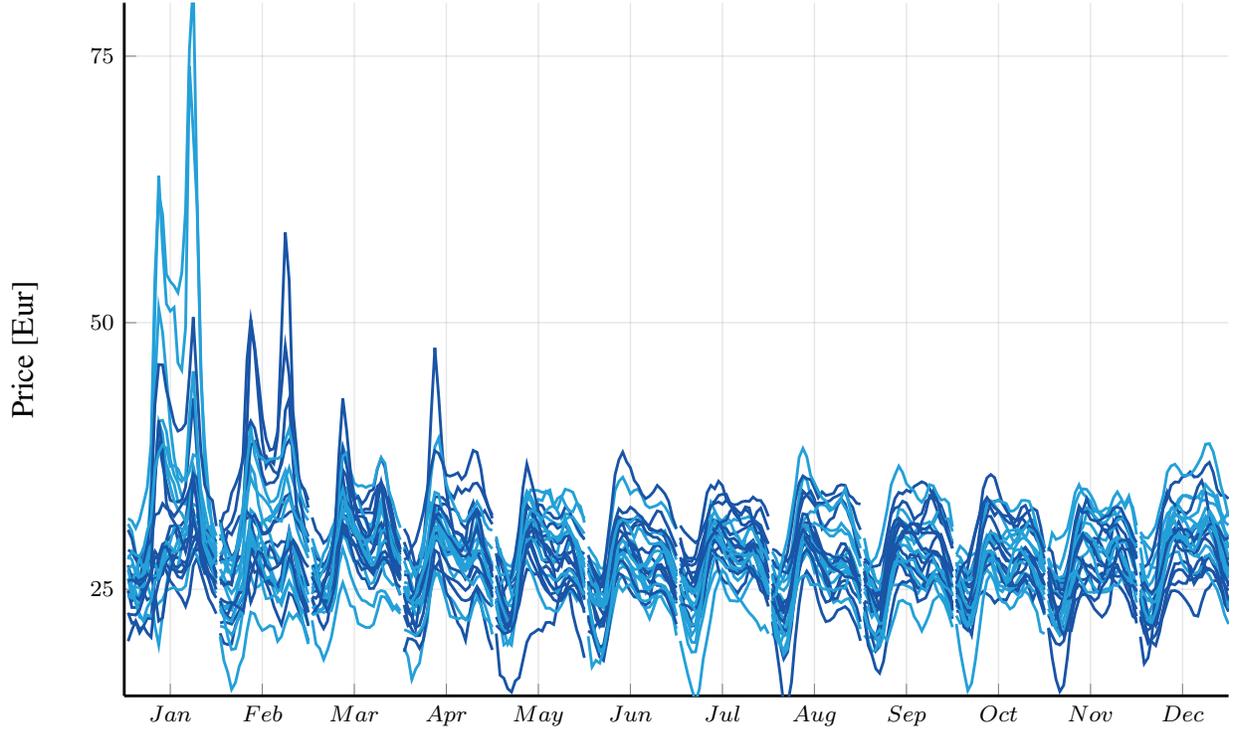}
  \caption{Daily electricity price curves predicted by the RNN forecaster in every month of the year.}
  \label{fig:seasonalforecast}
\end{figure}

We present the price forecaster and the inflow forecaster separately. Figure~\ref{fig:dayaheadforecast} shows historical price curves from January between the years 2013 and 2016, together with 200 price curves sampled using the RNN forecaster. The forecaster is set to predict curves in January. The neural network has been able to learn the overall shape of a typical electricity price curve. Price spikes typically occur in the morning when people wake up, and in the afternoon when people arrive home back from work. The network does not predict the extreme outliers of the past but does predict some daily curves with significant spikes. Next, we consider varying the month considered by the RNN forecaster. The result is shown in Figure~\ref{fig:seasonalforecast}. Each month, the shape of the daily price curve is consistent. However, the forecaster predicts larger prices in general in the winter period. This is consistent with the historical data, as electricity price in the Nordic region is higher in winter due to the heating demand. Similar results are obtained for the flow forecast.

\subsection{Day-ahead planning}
\label{sec:ch5-day-ahead-planning}

We have considered the computational performance of the algorithms in depth before~\cite{distlshaped}. Therefore, we refrain from that here and focus on the model output.

We solve the day-ahead problem for every month of the year by supplying an arbitrary date for each month to the RNN forecaster. By first solving some small sampled models, we note that the relative VSS is typically on the order of $\num{e-3}$. Therefore, we run the sequential SAA algorithm until the relative length of the confidence interval around the optimum is of order $\num{e-4}$ to a significance level of $95\%$. Consequently, the resulting VSS confidence intervals have a significance level of $90\%$. We initialize the algorithm at $n = 16$ number of samples and double the amount each iteration.

\begin{figure}
  \centering
  \input{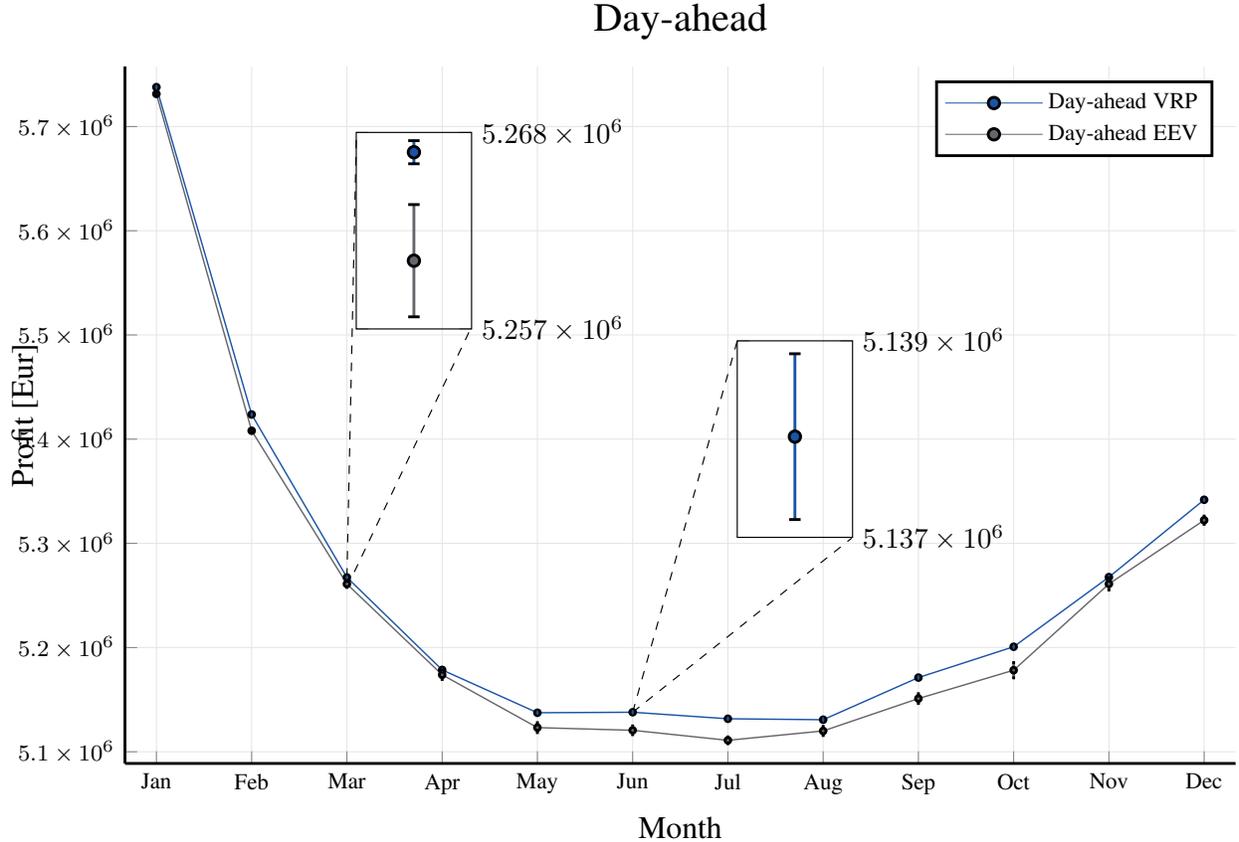}
  \caption{Seasonal variation of day-ahead VRP and EEV, including $95\%$ confidence intervals. Scale clarifications are shown for March and June. The VSS is statistically significant in all months except April and November.}
  \label{fig:seasonalvrp}
\end{figure}

\begin{figure}
  \centering
  \resizebox{\textwidth}{!}{
  \begin{tikzpicture}[]
    \begin{axis}[height = {101.6mm}, legend pos = {north west}, ylabel = {Profit [Eur]}, title = {Day-ahead}, xmin = {0.17000000000000004}, xmax = {11.83}, ymax = {31792.86786851776}, xlabel = {Month}, unbounded coords=jump,scaled x ticks = false,xlabel style = {font = {\fontsize{11 pt}{14.3 pt}\selectfont}, rotate = 0.0},xmajorgrids = true,xtick = {0.5,1.5,2.5,3.5,4.5,5.5,6.5,7.5,8.5,9.5,10.5,11.5},xticklabels = {Jan,Feb,Mar,Apr,May,Jun,Jul,Aug,Sep,Oct,Nov,Dec},xtick align = inside,xticklabel style = {font = {\fontsize{8 pt}{10.4 pt}\selectfont}, rotate = 0.0},x grid style = {color = kth-lightgray,
        line width = 0.25,
        solid},axis x line* = left,x axis line style = {line width = 1,
        solid},scaled y ticks = false,ylabel style = {font = {\fontsize{11 pt}{14.3 pt}\selectfont}, rotate = 0.0},ymajorgrids = true,ytick = {0.0,10000.0,20000.0,30000.0},yticklabels = {$0$,$1\times10^{4}$,$2\times10^{4}$,$3\times10^{4}$},ytick align = inside,yticklabel style = {font = {\fontsize{8 pt}{10.4 pt}\selectfont}, rotate = 0.0},y grid style = {color = kth-lightgray,
        line width = 0.25,
        solid},axis y line* = left,y axis line style = {line width = 1,
        solid},    xshift = 0.0mm,
      yshift = 0.0mm,
      title style = {font = {\fontsize{14 pt}{18.2 pt}\selectfont}, rotate = 0.0},legend style = {line width = 1,
        solid,font = {\fontsize{8 pt}{10.4 pt}\selectfont}},colorbar style={title=}, ymin = {-1396.1679265762214}, width = {152.4mm}]\addplot+[draw=none, color = kth-blue,
      line width = 0,
      solid,mark = *,
      mark size = 2.0,
      mark options = {
        color = black,
        fill = kth-blue,
        line width = 1,
        rotate = 0,
        solid
      }] coordinates {
        (0.5, 6486.015229139011)
        (1.5, 15584.260067903437)
        (2.5, 6349.431612970773)
        (3.5, 4681.823131513316)
        (4.5, 14268.114797750954)
        (5.5, 17334.704068879597)
        (6.5, 20670.731631594244)
        (7.5, 10773.313795404509)
        (8.5, 20145.898245697375)
        (9.5, 22511.117370738182)
        (10.5, 6544.446470947005)
        (11.5, 19602.380639248528)
      };
      \addlegendentry{Day-ahead VSS}
      \addplot+ [color = kth-blue,
      line width = 1,
      solid,mark = -,
      mark size = 2.0,
      mark options = {
        color = black,
        line width = 1,
        rotate = 0,
        solid
      },forget plot]coordinates {
        (0.5, 4226.420279428363)
        (0.5, 8745.61017884966)
      };
      \addplot+ [color = kth-blue,
      line width = 1,
      solid,mark = -,
      mark size = 2.0,
      mark options = {
        color = black,
        line width = 1,
        rotate = 0,
        solid
      },forget plot]coordinates {
        (1.5, 13494.003483271226)
        (1.5, 17674.516652535647)
      };
      \addplot+ [color = kth-blue,
      line width = 1,
      solid,mark = -,
      mark size = 2.0,
      mark options = {
        color = black,
        line width = 1,
        rotate = 0,
        solid
      },forget plot]coordinates {
        (2.5, 2388.1624825056642)
        (2.5, 10310.700743435882)
      };
      \addplot+ [color = kth-blue,
      line width = 1,
      solid,mark = -,
      mark size = 2.0,
      mark options = {
        color = black,
        line width = 1,
        rotate = 0,
        solid
      },forget plot]coordinates {
        (3.5, -456.85559275280684)
        (3.5, 9820.50185577944)
      };
      \addplot+ [color = kth-blue,
      line width = 1,
      solid,mark = -,
      mark size = 2.0,
      mark options = {
        color = black,
        line width = 1,
        rotate = 0,
        solid
      },forget plot]coordinates {
        (4.5, 8579.747557581402)
        (4.5, 19956.482037920505)
      };
      \addplot+ [color = kth-blue,
      line width = 1,
      solid,mark = -,
      mark size = 2.0,
      mark options = {
        color = black,
        line width = 1,
        rotate = 0,
        solid
      },forget plot]coordinates {
        (5.5, 11687.932158820331)
        (5.5, 22981.475978938863)
      };
      \addplot+ [color = kth-blue,
      line width = 1,
      solid,mark = -,
      mark size = 2.0,
      mark options = {
        color = black,
        line width = 1,
        rotate = 0,
        solid
      },forget plot]coordinates {
        (6.5, 16512.040094533935)
        (6.5, 24829.423168654554)
      };
      \addplot+ [color = kth-blue,
      line width = 1,
      solid,mark = -,
      mark size = 2.0,
      mark options = {
        color = black,
        line width = 1,
        rotate = 0,
        solid
      },forget plot]coordinates {
        (7.5, 5184.81263900362)
        (7.5, 16361.814951805398)
      };
      \addplot+ [color = kth-blue,
      line width = 1,
      solid,mark = -,
      mark size = 2.0,
      mark options = {
        color = black,
        line width = 1,
        rotate = 0,
        solid
      },forget plot]coordinates {
        (8.5, 14396.966961139813)
        (8.5, 25894.829530254938)
      };
      \addplot+ [color = kth-blue,
      line width = 1,
      solid,mark = -,
      mark size = 2.0,
      mark options = {
        color = black,
        line width = 1,
        rotate = 0,
        solid
      },forget plot]coordinates {
        (9.5, 14168.67920678202)
        (9.5, 30853.555534694344)
      };
      \addplot+ [color = kth-blue,
      line width = 1,
      solid,mark = -,
      mark size = 2.0,
      mark options = {
        color = black,
        line width = 1,
        rotate = 0,
        solid
      },forget plot]coordinates {
        (10.5, -413.77439679112285)
        (10.5, 13502.667338685133)
      };
      \addplot+ [color = kth-blue,
      line width = 1,
      solid,mark = -,
      mark size = 2.0,
      mark options = {
        color = black,
        line width = 1,
        rotate = 0,
        solid
      },forget plot]coordinates {
        (11.5, 14939.419726626948)
        (11.5, 24265.341551870108)
      };
    \end{axis}
  \end{tikzpicture}}
  \caption{Seasonal variation of day-ahead VSS, including $90\%$ confidence intervals. Note that the VSS is not statistically significant in April and November.}
  \label{fig:seasonalvss}
\end{figure}
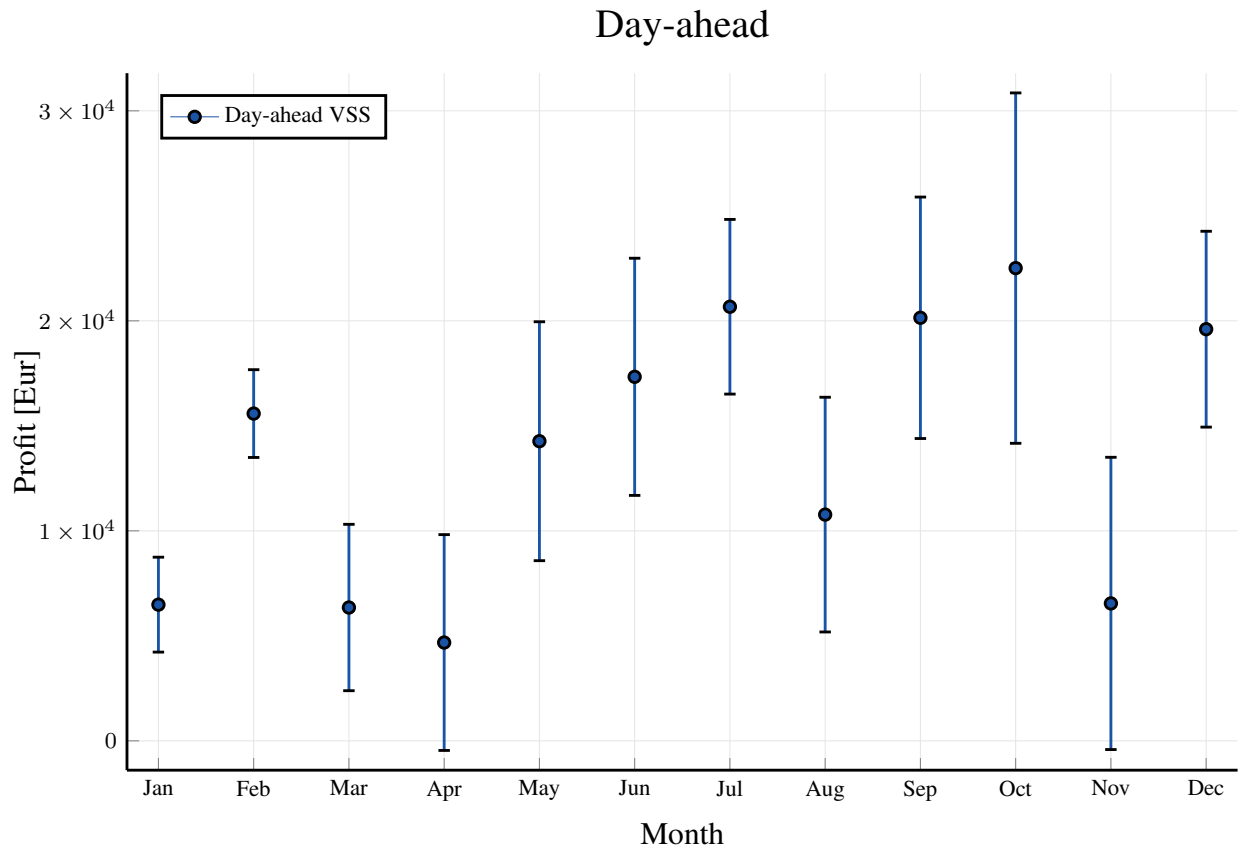

Every month, the SAA algorithm terminates at $n \approx 2000$. The resulting VRP and EEV confidence intervals are presented in Figure~\ref{fig:seasonalvrp}. The confidence intervals around the VRP values are tight. For example, the relative length of the VRP confidence interval in March is $\num{9.73e-5}$. A statistically significant VSS is ensured to the given confidence level as long as the confidence intervals around VRP and EEV do not overlap. It is no surprise that tight intervals are required here as there is an overlap in both April and November even at this low relative tolerance. In all other months, we can calculate a statistically significant VSS value. These seasonal VSS values are presented in more detail in Figure~\ref{fig:seasonalvss}. Notably, the relative VSS is small, ranging between $0.1\% - 0.4\%$. However, the value function includes the water evaluation, which increases the order of magnitude significantly. With respect to only the daily market profit, the relative VSS is about $1\%$. Besides, these constitute daily marginal profits. Hence, the VSS accumulates and could be considered more significant if the stochastic programming approach is employed over a longer time period.

\begin{figure}
  \centering
  \input{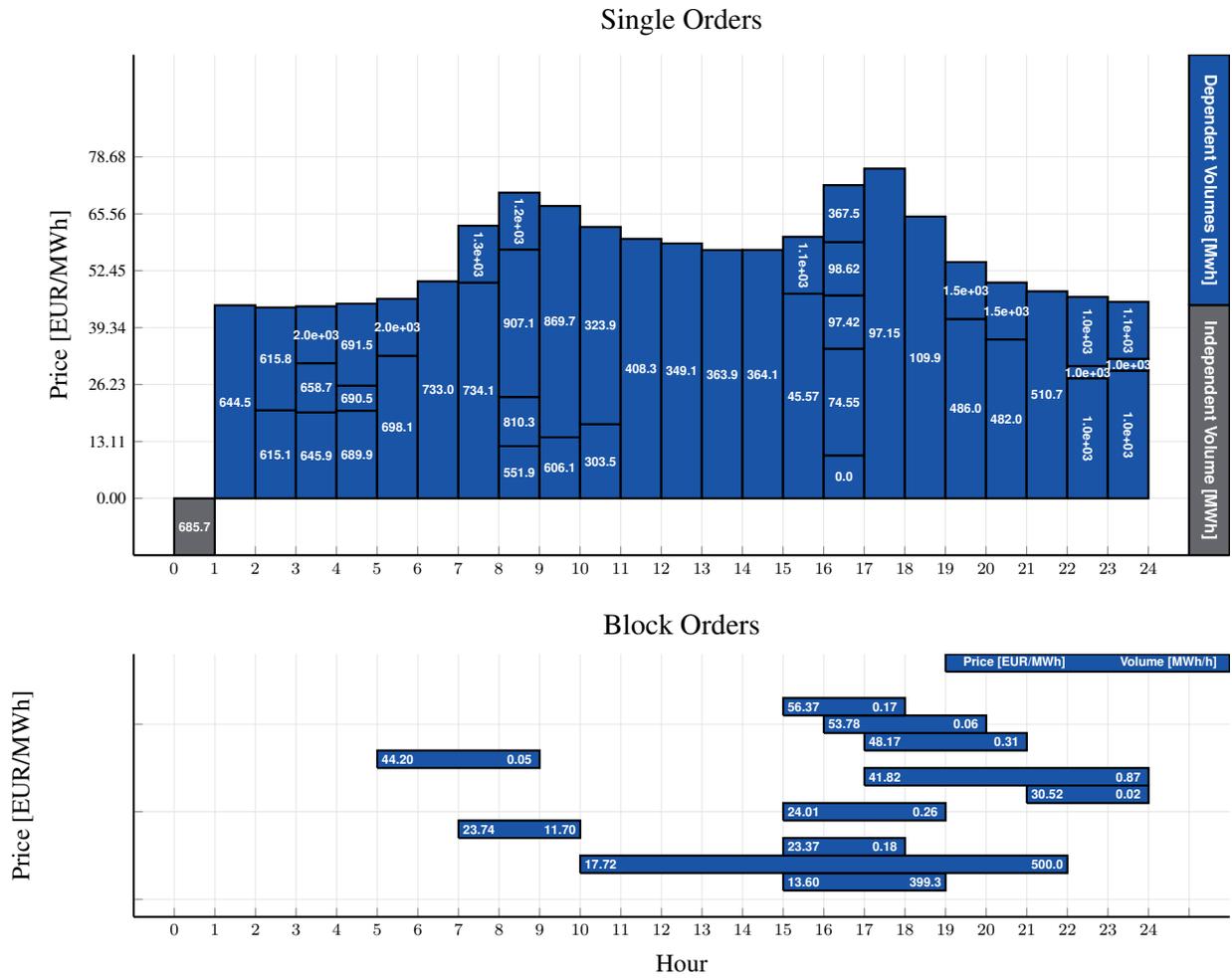}
  \caption{Optimal order strategy when solving a $2000$-scenario day-ahead SAA instance in January.}
  \label{fig:strategy}
\end{figure}

The optimal order strategy in January, when solving a $2000$-scenario day-ahead SAA instance, is shown in Fig.~\ref{fig:strategy}. The stochastic solution uses a large block order in the afternoon where a large mean price is expected. In comparison, the deterministic strategy obtained by solving the expected value problem is shown in Fig.~\ref{fig:evp_strategy}. The deterministic strategy mostly utilizes price-independent orders, which is less flexible than the stochastic solution.

\begin{figure}
  \centering
  \input{day-ahead_evp_strategy.tex}
  \caption{Optimal order strategy when solving the day-ahead expected value problem in January.}
  \label{fig:evp_strategy}
\end{figure}

\begin{figure}
  \centering
  \resizebox{\textwidth}{!}{
  \begin{tikzpicture}[]
    \begin{axis}[height = {101.6mm}, legend pos = {north west}, ylabel = {Profit [Eur]}, title = {Day-ahead}, xmin = {-30.0}, xmax = {1030.0}, ymax = {297792.1205850209}, xlabel = {Penalty percentage [\%]}, unbounded coords=jump,scaled x ticks = false,xlabel style = {font = {\fontsize{11 pt}{14.3 pt}\selectfont}, rotate = 0.0},xmajorgrids = true,xtick = {0.0,250.0,500.0,750.0,1000.0},xticklabels = {$0$,$250$,$500$,$750$,$1000$},xtick align = inside,xticklabel style = {font = {\fontsize{8 pt}{10.4 pt}\selectfont}, rotate = 0.0},x grid style = {color = kth-lightgray,
        line width = 0.25,
        solid},axis x line* = left,x axis line style = {line width = 1,
        solid},scaled y ticks = false,ylabel style = {font = {\fontsize{11 pt}{14.3 pt}\selectfont}, rotate = 0.0},ymajorgrids = true,ytick = {0.0,50000.0,100000.0,150000.0,200000.0,250000.0},yticklabels = {$0$,$5.0\times10^{4}$,$1.0\times10^{5}$,$1.5\times10^{5}$,$2.0\times10^{5}$,$2.5\times10^{5}$},ytick align = inside,yticklabel style = {font = {\fontsize{8 pt}{10.4 pt}\selectfont}, rotate = 0.0},y grid style = {color = kth-lightgray,
        line width = 0.25,
        solid},axis y line* = left,y axis line style = {line width = 1,
        solid},    xshift = 0.0mm,
      yshift = 0.0mm,
      title style = {font = {\fontsize{14 pt}{18.2 pt}\selectfont}, rotate = 0.0},legend style = {line width = 1,
        solid,font = {\fontsize{8 pt}{10.4 pt}\selectfont}},colorbar style={title=}, ymin = {-18859.953672402808}, width = {152.4mm}]\addplot+[draw=none, color = kth-blue,
      line width = 0,
      solid,mark = *,
      mark size = 2.0,
      mark options = {
        color = black,
        fill = kth-blue,
        line width = 1,
        rotate = 0,
        solid
      }] coordinates {
        (0.0, 2172.628359590657)
        (10.0, 3101.2972660232335)
        (50.0, 11208.487535539549)
        (100.0, 73831.77412320348)
        (500.0, 204726.50803605793)
        (1000.0, 269658.9830141496)
      };
      \addlegendentry{Day-ahead VSS}
      \addplot+ [color = kth-blue,
      line width = 1,
      solid,mark = -,
      mark size = 2.0,
      mark options = {
        color = black,
        line width = 1,
        rotate = 0,
        solid
      },forget plot]coordinates {
        (0.0, -7180.118704919703)
        (0.0, 11525.375424101017)
      };
      \addplot+ [color = kth-blue,
      line width = 1,
      solid,mark = -,
      mark size = 2.0,
      mark options = {
        color = black,
        line width = 1,
        rotate = 0,
        solid
      },forget plot]coordinates {
        (10.0, -9898.102514173836)
        (10.0, 16100.697046220303)
      };
      \addplot+ [color = kth-blue,
      line width = 1,
      solid,mark = -,
      mark size = 2.0,
      mark options = {
        color = black,
        line width = 1,
        rotate = 0,
        solid
      },forget plot]coordinates {
        (50.0, 1834.6670053666458)
        (50.0, 20582.308065712452)
      };
      \addplot+ [color = kth-blue,
      line width = 1,
      solid,mark = -,
      mark size = 2.0,
      mark options = {
        color = black,
        line width = 1,
        rotate = 0,
        solid
      },forget plot]coordinates {
        (100.0, 64878.6534786541)
        (100.0, 82784.89476775285)
      };
      \addplot+ [color = kth-blue,
      line width = 1,
      solid,mark = -,
      mark size = 2.0,
      mark options = {
        color = black,
        line width = 1,
        rotate = 0,
        solid
      },forget plot]coordinates {
        (500.0, 197175.1208734857)
        (500.0, 212277.89519863017)
      };
      \addplot+ [color = kth-blue,
      line width = 1,
      solid,mark = -,
      mark size = 2.0,
      mark options = {
        color = black,
        line width = 1,
        rotate = 0,
        solid
      },forget plot]coordinates {
        (1000.0, 250487.69660150725)
        (1000.0, 288830.2694267919)
      };
    \end{axis}

  \end{tikzpicture}

}
  \caption{Day-ahead VSS as a function of imbalance penalty, including $90\%$ confidence intervals calculated susing $1000$ scenarios. The penalty is given as a percentage of the market price. The VSS is not statistically significant for low penalty values.}
  \label{fig:penaltyvss}
\end{figure}
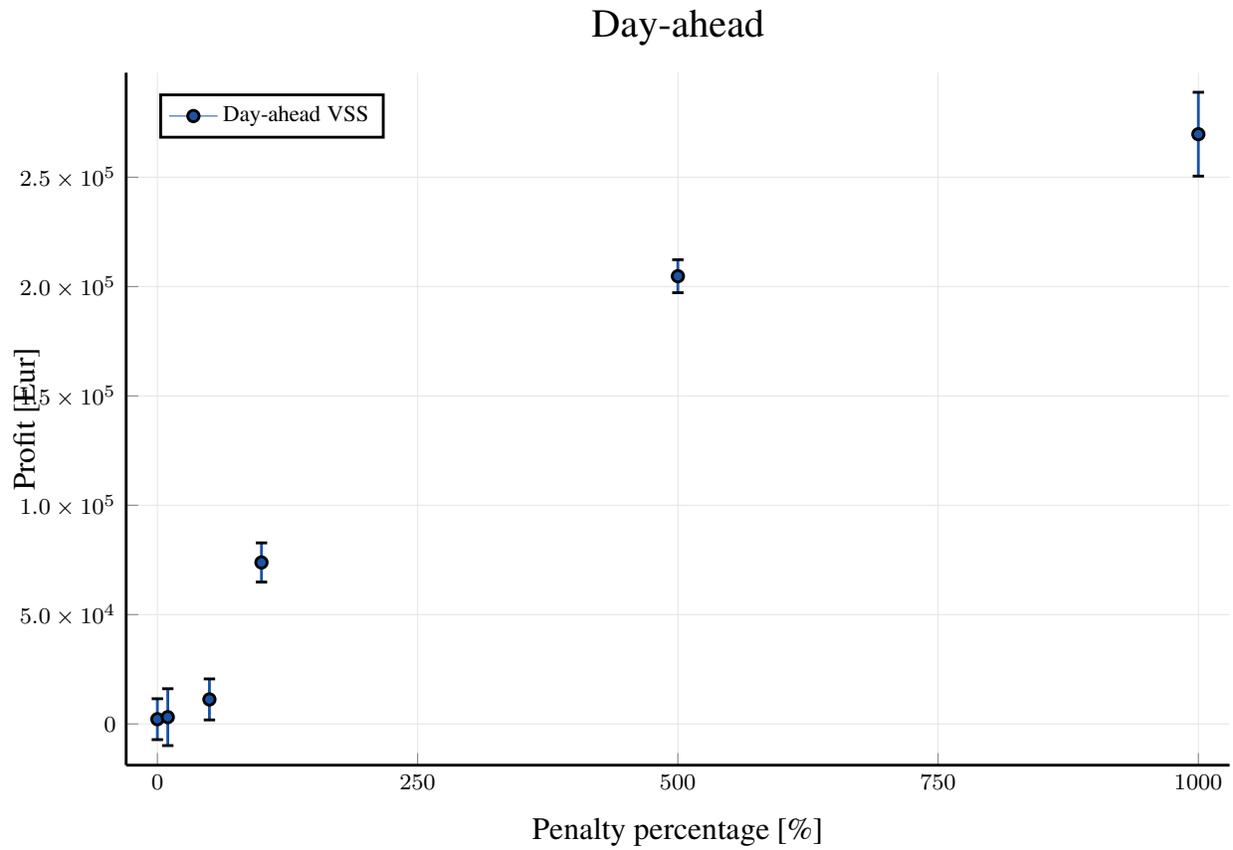

\begin{figure}
  \centering
  \input{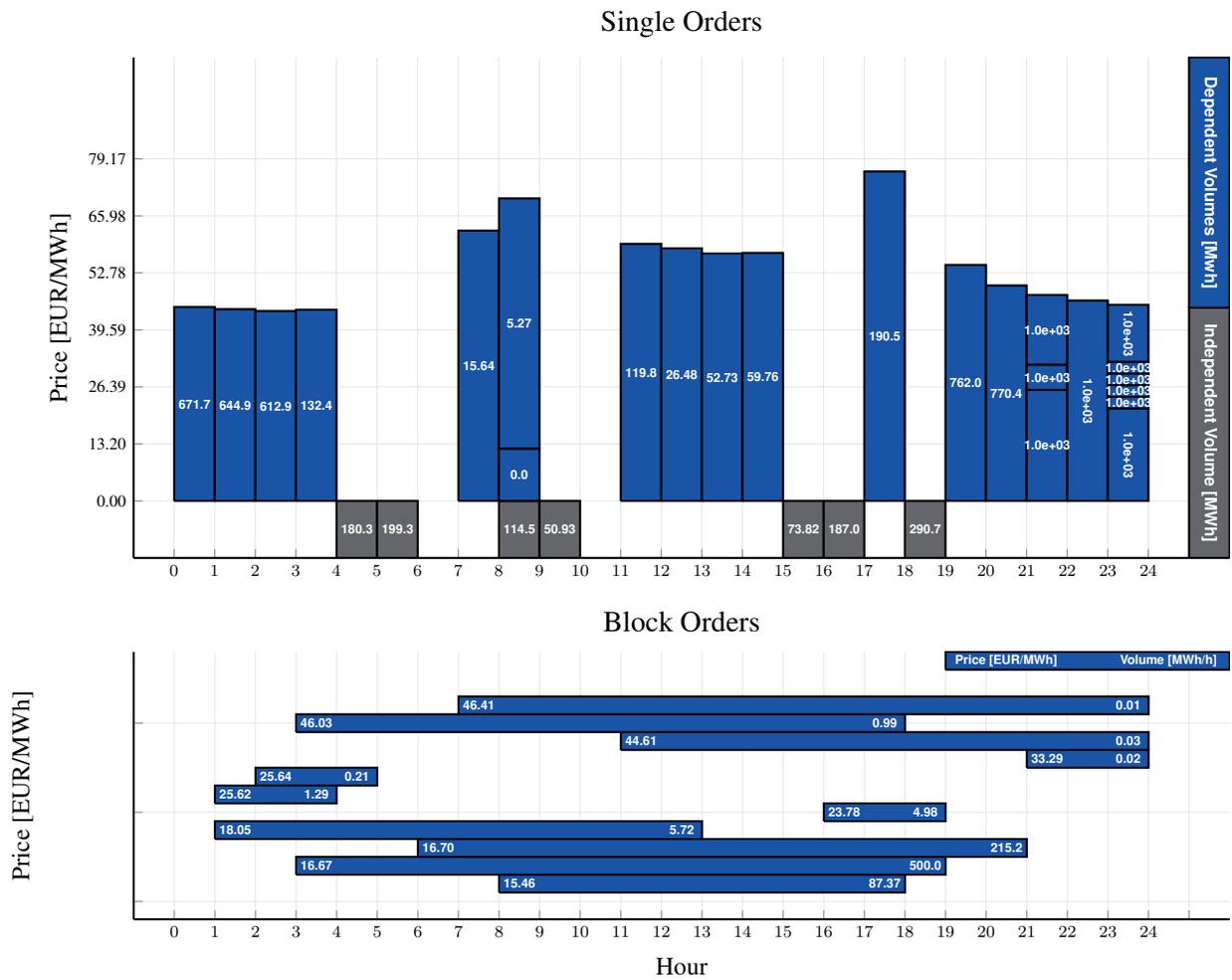}
  \caption{Optimal order strategy when solving a $1000$-scenario day-ahead SAA instance, with $1000\%$ imbalance penalty, in January.}
  \label{fig:penaltystrategy}
\end{figure}

The day-ahead VSS was shown to be sensitive to the imbalance penalty~\cite{prosumer}. We investigate this by letting the penalty be an increasing percentage of the market price in each scenario. An imbalance penalty of more than $100\%$ is not realistic, but can be interprested as a penalty formulation of a zero imbalance constraint. This models the case where no imbalance in the market commitments is desired, which could be beneficial to some producers. For larger penalty values, a statistically significant VSS can be achieved with fewer scenarios. Here, $1000$ scenarios are used to calculate the VSS values presented in Fig.~\ref{fig:penaltyvss}. It should be noted that the stochastic program is more complicated and takes about twice as long to solve with large imbalance penalties. The resulting order strategy is shown in Fig.~\ref{fig:penaltystrategy} and the deterministic counterpart is shown in Fig.~\ref{fig:evppenaltystrategy}. The stochastic solution is far more conservative to reduce the risk of imbalance. The deterministic strategy places larger orders and still makes heavy use of price-independent orders. The risk of imbalance is increased with this strategy, which is reflected in the VSS.

\begin{figure}
  \centering
  \input{day-ahead_penalty_evp_strategy.tex}
  \caption{Optimal order strategy when solving the day-ahead expected value problem, with $1000\%$ imbalance penalty, in January.}
  \label{fig:evppenaltystrategy}
\end{figure}

\section{Discussion and conclusion}
\label{sec:ch5-conclusion}

The results presented in this work are encouraging. However, the results are only as accurate as the water evaluation. The approach presented here for valuing the water is more involved than just assuming that all excess water can be sold at some conjectured future market price. However, it is probably not an accurate indicator of the water value. After scenario realization, the week-ahead problem we formulate assumes complete knowledge about the price curves and inflows for the upcoming week. A more precise formulation would involve a multi-stage stochastic program where the uncertain prices and inflows are learned sequentially. Another factor is the assumption that all power stations in the river are owned by the same fictitious hydropower producer. Skellefteälven is in reality operated by three separate companies and production schedules need to be coordinated. Hence, game-theoretic approaches are required for optimal planning. Finally, the RNN forecasters are trained on historical price and inflow data independently. In the Nordic region, hydropower constitutes more than half of the total electricity production. Consequently, the water value has a large impact on the electricity price. Thus, one should arguably consider training one single RNN forecaster on both datasets simultaneously, as electricity price and local inflow are most probably correlated.

In summary, a novel forecasting technique based on noise-driven recurrent neural networks has been presented. A key feature is that accurate forecasts can be generated without using lagged input sequences from historical data. We employ this forecasting technique to generate samples of day-ahead electricity prices and local reservoir inflows. We then formulate and solve a stochastic program for determining optimal day-ahead order strategies in the Nordic market, using the forecasters for scenario generation. The RNN forecasters predict accurate seasonal trends from historical data using only Gaussian signals and season indicators as input. The SAA algorithm yields tight confidence intervals around the stochastic solution. This allows us to derive statistically significant VSS of the stochastic order strategies for most months of the year.

\clearpage
\bibliographystyle{unsrt}
\bibliography{references}

\end{document}